\documentclass[11pt,a4paper]{article}

\usepackage{amsmath,amssymb, bbm}
\usepackage{color}
\newcommand{\R}{\mathbb{R}}
\newcommand{\N}{\mathbb{N}}
\newcommand{\Z}{\mathbb{Z}}
\newcommand{\A}{\mathbb{A}}
\newcommand{\B}{\mathbb{B}}
\def\cA{{\mathcal A}}
\def\cC{{\mathcal C}}
\def\cB{{\mathcal B}}
\def\cD{{\mathcal D}}
\def\cS{{\mathcal S}}
\def\cF{{\mathcal F}}
\def\cP{{\mathcal P}}
\def\cQ{{\mathcal Q}}
\newcommand{\ee}{\varepsilon}
\renewcommand{\aa}{\alpha}
\renewcommand{\div}{{\rm div}\,}

\def\d{\partial}

\def\ddq{\dot \Delta_q}

\def\tilde{\widetilde}
\def\hat{\widehat}
\newcommand{\D}{\Delta}
\newcommand{\DF}{\Delta_F}
\newcommand{\La}{\Lambda}
\newcommand{\n}{\nabla}
\newcommand{\G}{\Gamma}
\newcommand{\Ge}{G^e}
\newcommand{\Fe}{F^e}
\newcommand{\Om}{\Omega}
\newcommand{\Ome}{\Omega_\varepsilon}
\newcommand{\GL}{\Gamma_L}

\newcommand{\fr}{\frac{1}{r}}
\newcommand{\frb}{\frac{2}{r}}

\newcommand{\no}{\nu_0}

\newcommand{\ve}{v_\ee}
\newcommand{\Ue}{U_\ee}
\newcommand{\Uqg}{U_{\ee,QG}}
\newcommand{\Uosc}{U_{\ee, osc}}
\newcommand{\Phie}{\Phi_\ee}
\newcommand{\Thee}{\theta_\ee}

\newcommand{\Xtl}{X_{t,\lambda}^\ee}
\newcommand{\Xtol}{X_{\tau,\lambda}^\ee}
\newcommand{\Xtli}{X_{t,\lambda}^{\ee,i}}
\newcommand{\Xdo}{X^\ee(x,D) \Omega_\ee}
\newcommand{\M}{M_p}
\newcommand{\cu}{\underbar{c}}

\newtheorem{thm}{Theorem}
\newtheorem{lem}{Lemma}

\newtheorem{prop}{Proposition}
\newtheorem{defi}{Definition}

\newtheorem{rem}{Remark}

\title{Asymptotics and lower bound for the lifespan of solutions to the Primitive Equations}

\author{Fr\'ed\'eric Charve\footnote{Universit\'e Paris-Est Cr\'eteil, Laboratoire d'Analyse et de Math\'ematiques Appliqu\'ees (UMR 8050), 61 Avenue du G\'en\'eral de Gaulle, 94 010 Cr\'eteil Cedex (France). E-mail: frederic.charve@u-pec.fr}}

\date{}
\begin{document}

\maketitle

\begin{abstract} This article generalizes a previous work in which the author obtained a large lower bound for the lifespan of the solutions to the Primitive Equations, and proved convergence to the 3D quasi-geostrophic system for general and ill-prepared (possibly blowing-up) initial data that are regularization of vortex patches related to the potential velocity. These results were obtained for a very particular case when the kinematic viscosity $\nu$ is equal to the heat diffusivity $\nu'$, turning the diffusion operator into the classical Laplacian. Obtaining the same results without this assumption is much more difficult as it involves a non-local diffusion operator. The key to the main result is a family of a priori estimates for the 3D-QG system that we obtained in a companion paper.
\end{abstract}

\section{Introduction}

\subsection{Presentation of the physical models}

The Primitive Equations describe geophysical flows located in a large scale at the surface of the Earth (typically the atmosphere or an ocean) under the assumption that the vertical motion is much smaller than the horizontal one. Two important phenomena have to be considered in this case: the rotation of the Earth around its axis and the vertical stratification of the density induced by gravity.

When the motion is observed in a frame which is located at the surface, the rotation induces two additional terms in the equations: the Coriolis force and the centrifugal force. The latter is included in the pressure term and forms the geopotential $\Phie$. The former induces a vertical rigidity in the fluid as described by the Taylor-Proudman theorem: under a fast rotation the velocity of all particles located on the same vertical is horizontal and constant. The influence of the rotation on a fluid motion depends on the fluid time-space scale and is measured by the Rossby number $Ro$ which basically compares the frequency of the rotation to the characteristic time of the motion.

Gravity induces a horizontal rigidity to the fluid density: heavier masses lay under lighter ones. Internal motions in the fluid tend to alter this layered structure but gravity tends to restore it, and the importance of this rigidity is measured by the Froude number $Fr$, directly related to the Brunt-V\"ais\"al\"a frequency and the buoyancy.

The smaller are these numbers, the more important are these phenomena and in the present paper, we will consider the Primitive Equations in the whole space and for small $Ro$ and $Fr$ numbers under the same scale: $Ro=\ee$ and $Fr=\ee F$ with $F\in[0,1]$. In what follows we will call $\ee$ the Rossby number and $F$ the Froude number. The system then writes as follows:
\begin{equation}
\begin{cases}
\d_t \Ue +\ve\cdot \n \Ue -L \Ue +\frac{1}{\ee} \cA \Ue=\frac{1}{\ee} (-\n \Phie, 0),\\
\div \ve=0,\\
{\Ue}_{|t=0}=U_{0,\ee}.
\end{cases}
\label{PE}
\tag{$PE_\ee$}
\end{equation}
The unknows are $\Ue =(\ve, \Thee)=(\ve^1, \ve^2, \ve^3, \Thee)$ ($\ve$ denotes the velocity of the fluid and $\Thee$ the scalar potential temperature, strongly related to the density fluctuation) and $\Phie$ the geopotential. The diffusion operator $L$ is defined by
$$
L\Ue \overset{\mbox{def}}{=} (\nu \D \ve, \nu' \D \Thee),
$$
where $\nu, \nu'>0$ are the cinematic viscosity and the thermal diffusivity. The matrix $\cA$ is defined by
$$
\cA \overset{\mbox{def}}{=}\left(
\begin{array}{llll}
0 & -1 & 0 & 0\\
1 & 0 & 0 & 0\\
0 & 0 & 0 & F^{-1}\\
0 & 0 & -F^{-1} & 0
\end{array}
\right)
$$
We also assume that the sequence of initial data is convergent as $\ee$ goes to zero:
$$
U_{0,\ee} \underset{\ee \rightarrow 0}{\longrightarrow} U_0=(v_0, \theta_0).
$$
This system generalises the following rotating fluid system, which adds to the classical incompressible Navier-Stokes system the Coriolis force and the centrifugal force ($\ee$ is the Rossby number):
\begin{equation}
\begin{cases}
\d_t \ve+\ve\cdot \n \ve -\nu \D \ve +\frac{\ve\times e_3}{\ee}=-\n p_\ee,\\
\div \ve =0,\\
{\ve}_{|t=0}=v_0
\end{cases}
\label{RF}
\tag{$RF_\ee$}
\end{equation}
We refer to \cite{BMN5, Cushman, EmMa, EmMa2, FC5, Greenspan, Sadourny, Pedlosky, LiTeWa1, LiTeWa2} for a more precise presentation of the physical models, to \cite{BMN1, BMN4, CDGG, CDGG3, Dutrifoy2, IG1} concerning the rotating fluids system. We also mention \cite{DesGre, Dragos1, CDGG2, Marius1, Marius2, Marius3, VSN, FC5} for results with anisotropic (and possibly evanescent) viscosities and \cite{Yoneda} for the stationnary rotating fluid system.
\\

Dividing $\cA \Ue$ or $\ve\times e_3$ by $\ee$ imposes conditions for the limit as $\ee$ goes to $0$, these terms are said to be penalized. Compared to the classical Navier-Stokes system (NS), \eqref{PE} and \eqref{RF} introduce a penalized additional term involving a skew-symmetric matrix, so that for the cannonical scalar product we have $\cA \Ue\cdot \Ue=(\ve\times e_3)\cdot \ve=0$. Therefore any energy method (relying on energy estimates in $L^2$ or in $H^s/\dot{H}^s$) will not see these penalized terms and will work as for $(NS)$. Then the Leray and Fujita-Kato theorems are very easily adapted and provide global (unique in $2D$) weak solutions if $U_{0,\ee}\in L^2$ and local unique strong solutions if $U_{0,\ee} \in \dot{H}^\frac{1}{2}$ (global for small initial data). In the present paper we will consider ill-prepared initial data (with large and blowing up norms as $\ee$ goes to zero) for which these fundamental results cannot give us anymore information.
\\

As mentionned before, we are interested in the asymptotics, as the small parameter $\ee$ goes to zero (implying a competition between the vertical structure induced by the rotation and the horizontal structure induced by gravity). As explained in \cite{FC, FC2}, filtering the fast oscillations helps stabilizing the system: if the rotation and stratification are strong enough (that is when $\ee$ is small), then \eqref{PE} has global solutions and all that remains is the slow motion satisfying the following limit system which is called the 3D quasi-geostrophic system (it was used in the first half of 20th century for modeling and forecasting at mid-latitude the atmospheric and oceanic circulation, it is now less used except for slow motion climate systems) and can be written as the coupling of a transport-diffusion equation with a Biot-Savart-type law as follows:
\begin{equation}
\begin{cases}
\d_t \Om +v.\n \Om -\G \Om =0\\
U=(v,\theta)=(-\partial_2, \partial_1, 0, -F\partial_3) \DF^{-1} \Om,\\
\Om_{|t=0}=\Om_0,
\end{cases}
\label{QG1}\tag{QG}
\end{equation}
where the operator $\G$ is defined by:
$$
\G \overset{def}{=} \D \DF^{-1} (\nu \d_1^2 +\nu \d_2^2+ \nu' F^2 \d_3^2),
$$
with $\DF=\d_1^2 +\d_2^2 +F^2 \d_3^2$ and where we also have $\Om=\d_1 U^2 -\d_2 U^1 -F \d_3 U^4 =\d_1 v^2 -\d_2 v^1 -F \d_3 \theta$ and $\Om_0=\d_1 v_0^2 -\d_2 v_0^1 -F \d_3 \theta_0$ where $U_0=(v_0, \theta_0)$ is the limit as $\ee$ goes to zero of the initial data $U_{0,\ee}$.
\begin{rem}
\sl{Except in the particular cases $\nu=\nu'$ (where $\G=\nu \D$, see \cite{FC3} and \cite{Dutrifoy2} for example) or $F=1$ (where $\G= \nu \d_1^2 +\nu \d_2^2+ \nu' \d_3^2$, see \cite{Chemin2}) the operator $\G$ is a non-local diffusion operator of order 2 and we refer to \cite{FCestimLp} for a detailed study.
}
\end{rem}
We emphasize that two cases have to be considered with respect to $F$: the non-dispersive case $F=1$, and the dispersive case $F\neq 1$ (this denomination is explained after Proposition \ref{estimvp}). In these settings is proved the convergence of the solutions of \eqref{PE} to the solution of \eqref{QG1} as the small parameter $\ee$ goes to zero. If $F=1$, in \cite{Chemin2} the result for regular well-prepared initial data is proved provided that $\nu$ and $\nu'$ are very close, in \cite{Dragos4} the convergence is proved in the inviscid case. In the case $F\neq 1$ another approach initiated by \cite{CDGG} consists in using dispersive and Strichartz estimates in order to filter fast oscillations which leads to a stabilization of the system and make it tend to the 2D Navier-Stokes system (in the case of the rotating fluid system, this is the Taylor-Proudman theorem) or to the slow 3D quasi-geostrophic system (in the case of the Primitive Equations). We obtained this convergence for various ill-prepared initial data (possibly blowing up) or anisotropic viscosities in \cite{FC, FC2, FC3, FC4, FC5}. As expected by physicists, we obtained that the difference of the solutions of the Primitive Equations and the quasi-geostrophic system is of the size of the Rossby number. We also refer to \cite{Dutrifoy2} for the inviscid case. We adress the reader to Remark \ref{illprepared} for the notion of  well or ill-prepared initial data.
\\

In \cite{FC} ($F\neq 1$) we first obtain formally the quasi-geostrophic system, and then, guided by its form we introduce the following decomposition: we first define the potential vorticity of any 4-dimensional vector field $U=(v, \theta)$ denoted by $\Om(U)$:
$$
\Om(U)\overset{def}{=} \d_1 v^2 -\d_2 v^1 -F \d_3 \theta.
$$
From this we define the quasi-geostrophic and oscillating parts of $U$:
\begin{equation}
U_{QG}=\cQ (U) \overset{def}{=} \left(
\begin{array}{r}
-\d_2\\
\d_1\\
0\\
-F\d_3
\end{array}
\right) \DF^{-1} \Om (U), \quad \mbox{and} \quad U_{osc}=\cP (U) \overset{def}{=} U-U_{QG}.
\end{equation}
As emphasized in \cite{FC,FC5} this is an orthogonal decomposition of 4-dimensional vectorfields (similar to the Leray orthogonal decomposition of any 3-dimensional vectorfield into its divergence-free and its gradient parts) and $\cQ$ and $\cP$ are the associated orthogonal projectors on the quasi-geostrophic or oscillating fields and satisfy (see \cite{Chemin2, FC, FC2}):
\begin{prop}
\sl{With the same notations, for any function $U=(v, \theta)$ we have the following properties:
\begin{enumerate}
\item $\cP$ and $\cQ$ are pseudo-differential operators of order 0,
\item For any $s\in\R$, $(\cP(U)|\cQ(U))_{\dot{H}^s} =(\cA U|\cP(U))_{\dot{H}^s}=0$,
\item The same is true for nonhomogeneous Sobolev spaces,
\item $\cP(U)=U \Longleftrightarrow \cQ(U)=0\Longleftrightarrow \Om(U)=0$,
\item $\cQ(U)=U \Longleftrightarrow \cP(U)=0\Longleftrightarrow$ there exists a scalar function $\Phi$ such that $U=(-\d_2,\d_1,0,-F\d_3) \Phi$. Such a vectorfield is said to be quasi-geostrophic and is divergence-free.
\item If $U=(v, \theta)$ is a quasi-geostrophic vectorfield, then $v\cdot \n \Om(U)=\Om(v\cdot \n U)$.
\item If $U$ is a quasi-geostrophic vectorfield, then $\G U=\cQ L U$.
\end{enumerate}
\label{propdecomposcqg}
}
\end{prop}
Thanks to this, we can rewrite the QG system into the following form:
\begin{equation}
\begin{cases}
\d_t U +\mathcal{Q} (v.\n U) -\G U=0,\\
U=\mathcal{Q} U,\\
U_{|t=0}=U_0.
\label{QG2}\tag{QG}
\end{cases}
\end{equation}
Going back to the Primitive Equations, if we denote by $\Ome=\Om(\Ue)$, $\Uqg=\cQ(\Ue)$ and $\Uosc=\cP(\Ue)$, they satisfy the following systems (see \cite{FC} for details):
\begin{equation}
\d_t \Ome +\ve\cdot \n \Ome -\G \Ome= (\nu-\nu')F\D \d_3 \theta_{\ee, osc}+q_\ee,
\label{systomega}
\end{equation}
where $q_\ee$ is defined by
\begin{multline}
q_\ee =q(\Uosc, \Ue) =\d_3 v_{\ee, osc}^3(\d_1 \ve^2-\d_2 \ve^1) -\d_1 v_{\ee, osc}^3 \d_3 \ve^2 +\d_2 v_{\ee, osc}^3 \d_3 \ve^1\\
+F\d_3 v_{\ee, QG}\cdot \n \theta_{\ee, osc} +F\d_3 v_{\ee, osc}\cdot \n \Thee,
\end{multline}
and
\begin{multline}
\d_t \Uosc-(L-\frac{1}{\ee}\mathbb{P} \cA) \Uosc = -\mathbb{P}(\ve\cdot \n \Ue) -\left(
\begin{array}{c}
-\d_2\\
\d_1\\
0\\
-F\d_3
\end{array}
\right) \DF^{-1} \bigg( -\ve\cdot \n \Ome +q(\Uosc, \Ue)\bigg)\\
+(\nu-\nu')F \D \D_F^{-1} \d_3\left(
\begin{array}{c}
\d_2 \theta_\ee\\
-\d_1 \theta_\ee\\
0\\
\d_1 v_\ee^2-\d_2 v_\ee^1
\end{array}
\right) \overset{def}{=} F_1 +F_2 +F_3 +F_4.
\label{systosc}
\end{multline}

\begin{rem}
\sl{For more simplicity and without any loss of generality, we will write in what follows
$$
q_\ee=\n \Uosc \cdot \n \Ue.
$$
}
\end{rem}

\begin{rem}
\sl{It is a natural question to know, when the initial data $U_{0,\ee}$ has a zero oscillating part (in other words, it is a purely quasi-geostrophic initial data), if the solution has the same property. From System \eqref{systosc} we can see that it may not be the case, but we can prove that the oscillating part is small in appropriate norms (see for example \cite{FC2}). We refer to \cite{Chemin2, Dragos4} where the results are given for initial data with small initial oscillating part (and vanishing as $\ee$ goes to zero). Such initial data are said to be \textit{well-prepared}, in contrast to the \textit{ill-prepared} case where the initial oscillating part can not only be large but also even blow-up as $\ee$ goes to zero as studied in the present paper. In this case it is essential to take advantage of the dispersive aspects of the system in order to control the oscillating part of the solution.
}
\label{illprepared}
\end{rem}

\subsection{Statement of the main results}

The aim of the present article is to use the estimates from the companion paper \cite{FCestimLp} and generalize the work from \cite{FC3} in the case $\nu \neq \nu'$. We refer to the appendix for general definitions about the Littlewood-Paley dyadic decomposition and the vortex patches formalism.

\begin{thm}
\sl{Assume that $U_{0, QG} \in L^2(\R^3)$ is a quasi-geostrophic vectorfield whose potential vorticity $\Om(U_{0, QG})=\Om_0 \in L^2(\R^3) \cap L^\infty(\R^3)$ is also a $C^s$-vortex patch with $s\in]0,1[$. Define $U_{0, \ee, QG}$ as the following regularization of $U_{0, QG}$:
$$
U_{0, \ee, QG} =\chi(\ee^\beta |D|) U_{0, QG}= \ee^{-3\beta} h(\ee^{-\beta}.)* U_{0, QG},
$$
where $\beta>0$, $\chi$ is a smooth cut-off function (see Section \ref{appendiceLP}) and $h=\cF^{-1}\left(\chi(|.|)\right)$. Let $(U_{0, \ee, osc})_{\ee>0}$ be a family of regular oscillating vectorfields (i.-e. with zero potential vorticity). Assume there exists $C_0>0$ such that for all $\ee>0$ the family of initial data $U_{0,\ee} =U_{0, \ee, QG}+U_{0, \ee, osc}$ satisfies:
\begin{equation}
\begin{cases}
\vspace{1mm}
\|U_{0, \ee}\|_{\dot{H}^1} \leq C_0,\\
\|U_{0, \ee}\|_{H^6} \leq C_0 \ee^{-5\beta}.
\end{cases}
\label{hypinit}
\end{equation}
Then there exists constants $C,\omega, \gamma, \ee_0>0$ ($\omega=\omega(F,s,\nu,\nu')$) such that if 
$$
\beta +C\gamma \leq \omega,
$$
then for all $\ee\leq \ee_0$, the lifespan $T_\ee^*$ of the solution $U_\ee$ satisfies: $T_\ee^* \geq T_\ee^\gamma \overset{def}{=}\gamma \ln(\ln|\ln\ee|)$ and for all $t\leq T_\ee^\gamma$ we have:
$$
\begin{cases}
\forall k\leq3,\quad \|\n^k \Uosc\|_{L_t^4 L^\infty} \leq \ee^{15 \omega},\\
V_\ee(t) \overset{def}{=} \displaystyle{\int_0^t \|\n \Ue (t')\|_{L^\infty} dt' \leq 2\gamma |\ln \ee|}.
\end{cases}
$$
Moreover we have local convergence: for all $T>0$, $\Uqg$ converges in $L^\infty([0,T], L^2)$ to the unique global (lipschitz) solution of the 3D-quasi-geostrophic system with initial data $U_{0, QG}$ (which is in $H^1(\R^3)$).
}
\label{thprincipal}
\end{thm}

\begin{rem}
\sl{We refer to \cite{FC2} for the fact that System $(QG)$ has a unique global solution if the initial data belongs to $H^1$.
}
\end{rem}

\begin{rem}
\sl{More precisely, to give an idea, we can choose $\omega=10^{-4}$.
}
\end{rem}

\begin{rem}
\sl{
We point out that the assumptions on the initial data imply that
$$
\begin{cases}
\vspace{1mm}
\|U_{0, \ee, QG} -U_{0, QG} \|_{\dot{H}^\sigma} \leq C \ee^{\beta(1-\sigma)} \|U_{0, QG} \|_{\dot{H}^1} \quad \mbox{if} \quad \sigma \in[0,1],\\
\displaystyle{\|U_{0, \ee, QG}\|_{\dot{H}^\sigma} \leq \frac{C}{\ee^{\beta(\sigma-1)}} \|U_{0, QG} \|_{\dot{H}^1} \quad \mbox{if} \quad \sigma>1.}
\end{cases}
$$
And concerning the initial oscillating part, it covers for example the case of the regularization of an oscillating function $U_{0, osc} \in \dot{H}^1$ which is not in $L^2$:
$$
U_{0, \ee, osc}=\psi(\ee^\beta |D|)U_{0, osc},
$$
where the function $\psi$ is supported in an annulus centered at $0$ and equal to $1$ in a smaller annulus.
}
\end{rem}
\begin{rem}
\sl{In \cite{FC3} we made a technical (and physically irrelevant) assumption: $\nu =\nu'$. In reality, and as suggested by the Prandtl number (which  can be defined as $\nu/\nu'$ and can take values far from $1$) there is no reason for the kinematic viscosity and the thermal diffusivity to be equal or even close (see also \cite{Chemin2}). We want, in the present article, to cover this general case $\nu \neq \nu'$, for which every step from \cite{FC3} will be much more tedious:
\begin{itemize}
\item First, the operator $\G$ is now a non-local operator. Dealing with products, commutators (and also commutator with a Lagrangian change of variable in \cite{FCestimLp}) will be much more difficult compared to the case of the classical Laplacian. A refined study of this operator is done in \cite{FCestimLp} in order to obtain the analogous of the a priori estimates used in \cite{FC3}.
\item When studying the persistence of the tangential (or stratified) regularity, as in \cite{TH1, FC3} we will have to cope with the commutator $[\Xtl(x,D),\G] \Ome$ whose study (when $\G$ is non-local) is much more difficult than in the case $\G=\nu \D$, and will require for example  more flexible estimates for the term $\G(fg)-f\G g-g\G f$ than what we needed in \cite{FCestimLp}.
\item Compared to \cite{FC3}, System \eqref{systomega}, satisfied by the potential vorticity $\Ome$, shows an additional external force term, involving derivatives of order 3, $(\nu-\nu')F\D \d_3 \theta_{\ee, osc}$ (which vanished thanks to the assumptions in \cite{FC3}). Because of this the general $L^\infty$-estimates for System \eqref{QG1} from \cite{FCestimLp} need to be modified, and we need a completely different approach for the $L^2$-estimates (we recall that the Leray estimates are useless as the $L^2$-norm of the initial data blows-up when $\ee$ goes to zero). This approach will strongly rely on the particular structure of the quasi-geostrophic decomposition (i.-e. the properties of operators $\cP$ and $\cQ$).
\item Compared to \cite{FC3}, to treat the general case $\nu \neq \nu'$, we need an additionnal assumption: $\|U_{0, \ee, osc}\|_{\dot{H}^1} \leq C_0$. This is once more due to the first term in the right-hand side of system \eqref{systomega}. Moreover, we emphasize that from the assumptions of the main result, the best we can hope for $L^2$ or $\dot{H}^\frac{1}{2}$ estimates are negative powers of $\ee$.
\item The $L^\infty$ estimates of the potentiel vorticity will introduce in the computations a multiplicative exponential factor $D^t$ that will degrade the size of the lower bound for the lifespan: in \cite{FC3} we obtained $T_\ee^*=\gamma \ln|\ln \ee|$. The consequence is that the frequency truncations of size $(-\ln \ee)^\delta$ from \cite{FC3} (when obtaining dispersive and Strichartz estimates) are useless in the present case and will have to be improved to the size $\ee^{-\delta}$.
\end{itemize}
}
\end{rem}

\section{Scheme of the proof and structure of the article}

As explained before, thanks to the skew-symmetry of matrix $\cA$, any computation involving $L^2$ or Sobolev inner-products will be the same as for the Navier-Stokes system ($\cA U\cdot U=0$). So given the regularity of the initial data (even if some norms can blow up in $\ee$), we can adapt the Leray and Fujita-Kato theorems as well as the classical weak-strong uniqueness results: as $U_{0,\ee}\in \dot{H}^\frac{1}{2}(\R^3)$, $\Ue$ is the unique strong solution of System \eqref{PE} defined on $[0,T]$ for all $0<T<T_\ee^*$. In addition, if the lifespan $T_\ee^*$ is finite then we have:
$$
\int_0^{T_\ee^*} \|\n \Ue(\tau)\|_{\dot{H}^\frac{1}{2}(\R^3)}^2 d\tau=\infty.
$$
Moreover, as $U_{0,\ee}\in L^2$, $\Ue$ coincides on $[0, T_\ee^*[$ with any global weak solution of System \eqref{PE} and we have the Leray estimates: for all $t\geq 0$,
\begin{equation}
\|\Ue(t)\|_{L^2}^2 + 2\nu_0\int_0^{t} \|\n \Ue (\tau)\|_{L^2}^2 d\tau \leq \|U_{0,\ee}\|_{L^2}^2.
\label{estimLeray}
\end{equation}
From \eqref{hypinit}, the only assumption on the initial regularity of $\Ue$ provides that $\|U_{0,\ee}\|_{\dot{H}^\frac{1}{2}(\R^3)} \leq C \ee^{-5\beta}$ which is useless if we want a lower bound for the Fujita-Kato lifespan without any other information. In this paper we will use the quasi-geostrophic structure to prove that the lifespan is bounded from below by some large time $T_\ee^\gamma = \gamma \ln(\ln|\ln \ee|)$. Moreover we will extensively use the following a priori estimates (we refer to \cite{Dutrifoy1, FC3}):
\begin{lem}
\sl{Let $s>1$. There exists a constant $C_s$ such that for all $t\in [0, T_\ee^*[$ we have:
\begin{equation}
\|\Ue(t)\|_{H^s} \leq \|\Ue (0)\|_{H^s} e^{C_s \int_0^t \|\n \Ue (\tau)\|_{L^\infty} d\tau},
\label{estapriori}
\end{equation}
}
\end{lem}
which requires us to be able to estimate $\|\n \Ue\|_{L_t^1 L^\infty}$. We will easily show in the following that $\|\n \Uosc\|_{L_t^1 L^\infty}$ is bounded and small. Showing that $\|\n \Uqg\|_{L_t^1 L^\infty}$ is bounded will require much more work, and as in \cite{TH1, Dutrifoy1, FC3} we will need the following logarithmic estimates, involving striated regularity from the vortex patches formalism (we refer to \cite{Dutrifoy1} for the proof) :
\begin{lem}
\sl{There exists a constant $C>0$ only depending on $s\in ]0,1[$ such that, for any quasi-geostrophic vector field $U\in L^2(\R^3)$ whose potential vorticity $\Om=\Om (U)\in L^2(\R^3) \cap C^s(X)$ for a fixed admissible family $X$ of $C^s$-vectorfields, $U$ is Lipschitzian and we have~:
\begin{equation}
\|U\|_{Lip}= \|\nabla U\|_{L^{\infty}} \leq C\Big(\|\Omega\|_{L^2}+ \|\Omega\|_{L^{\infty}} \log
\big(e+ \frac{\|\Omega\|_{C^s(X)}}{\|\Omega\|_{L^{\infty}}} \big) \Big).
\label{estimlog}
\end{equation}
}
\end{lem} 
The bootstrap argument is pretty simple: let us define $K_\ee\overset{def}{=} 2\gamma |\ln \ee|$, $T_\ee^\gamma \overset{def}{=} \gamma \ln(\ln|\ln \ee|)$ and the times $T_\ee> T_\ee'>0$ by
\begin{equation}
\begin{cases}
\vspace{1mm}
\displaystyle{T_\ee =\sup \{t\in [0,T_\ee^*[, \quad V_\ee(t)=\int_0^t \|\n \Ue(\tau)\|_{L^\infty} d\tau \leq K_\ee\},}\\
T_\ee' =\sup \{t\in [0,T_\ee^*[, \quad V_\ee(t) \leq K_\ee/2\}.
\end{cases}
\label{bootstrap}
\end{equation}
We will show in this article that for a well chosen $\gamma$ and if $\ee>0$ is small enough, for all $t\leq \min(T_\ee, T_\ee^\gamma)$, 
\begin{equation}
V_\ee(t) \leq K_\ee/2,
\label{bootestim}
\end{equation}
then, due to the previous definitions, $t \leq T_\ee'$, which immediately implies that $\min(T_\ee, T_\ee^\gamma) \leq T_\ee' < T_\ee$ and therefore $T_\ee^\gamma <T_\ee < T_\ee^*$ which proves the theorem.
\begin{rem}
\sl{Either $T_\ee^*<\infty$ and the integral goes to infinity when $t\rightarrow T_\ee^*$ and we have $0<T_\ee'< T_\ee$, or $T_\ee^*=\infty$ and then $T_\ee^\gamma < T_\ee^*$ is also true (in this case we may be have $T_\ee'= T_\ee$).}
\end{rem}

All the difficulty then lies in proving \eqref{bootestim} and the article is structured as follows, we will first show in Section $3$ that the oscillating part is small and goes to zero. As explained before, we will need to use logarithmic estimates \eqref{estimlog} related to the striated regularity which requires us to bound the $L^2$, $L^\infty$ and $C^s$ norms of the potential vorticity. The first one is obtained using the quasigeostrophic structure : in Section $4$ we get estimates in $\dot{H}^1$ for $\Ue$. The rest is dealt in Section $5$. As in \cite{FCestimLp} the major difficulty in this paper comes from the non-local operator $\G$ (we recall that in \cite{FC3} we only adressed the particular case $\nu=\nu'$ where $\G$ reduces to $\nu \D$). In this paper we will extensively use the apriori and smoothing estimates obtained in \cite{FCestimLp}. The last section is an appendix devoted to the Littlewood-Paley decomposition, followed by a quick presentation of the vortex patches formalism, and additional properties for the non-local operator $\G$.

\section{Estimates for the oscillating part}

This section is devoted to a careful study of the oscillating part. As in \cite{FC} to \cite{FC5} the fact that it goes to zero is essential in the study of the asymptotics of the quasi-geostrophic part and the convergence rates. The aim of this section is to prove that there exists $\delta>0$ such that if $\ee$ is small enough, for all $k\leq 3$,
$$
\||D|^k\Uosc\|_{L_T^4 L^\infty} \leq C_F \ee^{\delta}.
$$
\begin{rem}
\sl{In this article we will only need three derivatives, but it is true for any $k$.}
\end{rem}
We refer to Proposition \ref{estimUosc} for the precise statement of the result. As $\Uosc$ satisfies \eqref{systosc}, we will consider the following system:
\begin{equation}
\begin{cases}
\d_t f-(L-\frac{1}{\ee} \mathbb{P} \cA) f=\Fe,\\
f_{|t=0}=f_0.
\end{cases}
\label{systdisp}
\end{equation}
If we apply the Fourier transform, the equation becomes (see \cite{FC} for precisions):
$$
\d_t \hat{f}- \mathbb{B}(\xi, \ee)\hat{f}=\hat{\Fe},
$$
where
$$\mathbb{B}(\xi, \ee)= \hat{L-\frac{1}{\ee} \mathbb{P} \cA} =\left(
\begin{array}{cccc}
\displaystyle{-\nu|\xi|^2+\frac{\xi_1\xi_2}{\ee
  |\xi|^2}} & \displaystyle{\frac{\xi_2^2+\xi_3^2}{\ee
  |\xi|^2}} & 0 & \displaystyle{\frac{\xi_1\xi_3}{\ee F |\xi|^2}}\\
\displaystyle{-\frac{\xi_1^2+\xi_3^2}{\ee
  |\xi|^2}} & \displaystyle{-\nu|\xi|^2-\frac{\xi_1\xi_2}{\ee
  |\xi|^2}} & 0 & \displaystyle{\frac{\xi_2\xi_3}{\ee F |\xi|^2}}\\
\displaystyle{\frac{\xi_2\xi_3}{\ee |\xi|^2}} &
\displaystyle{-\frac{\xi_1\xi_3}{\ee
  |\xi|^2}} & \displaystyle{-\nu |\xi|^2} & \displaystyle{-\frac{\xi_1^2+\xi_2^2}{\ee F
  |\xi|^2}}\\
0 & 0 & \displaystyle{\frac{1}{\ee F}} &
\displaystyle{-\nu'|\xi|^2}
\end{array}
\right).
$$
For $0<r<R$ we will denote by $\cC_{r,R}$ the following set:
$$
\cC_{r,R} =\{\xi \in \R^3, \quad |\xi|\leq R \mbox{ and } |\xi_3|\geq r\}.
$$
We also introduce the following frequency truncation operator on $\cC_{r,R}$:
$$
\cP_{r,R}=\chi (\frac{|D|}{R})\big(1-\chi (\frac{|D_3|}{r})\big),
$$
where $\chi$ is the smooth cut-off function introduced before and ($\mathcal{F}^{-1}$ is the inverse Fourier transform):
$$
\chi (\frac{|D|}{R}) f=\mathcal{F}^{-1} \Big(\chi(\frac{|\xi|}{R}) \hat{f}(\xi)\Big) \quad \mbox{and} \quad \chi (\frac{|D_3|}{r}) f=\mathcal{F}^{-1} \Big(\chi(\frac{|\xi_3|}{r}) \hat{f}(\xi)\Big),
$$
and the following derivation operator:
$$
|D|^s f =\mathcal{F}^{-1} (|\xi|^s \hat{f}(\xi)).
$$
In what follows we will use it for particular radii $r_\ee=\ee^m$ and $R_\ee =\ee^{-M}$, where $m$ and $M$ will be precised later. Let us end this section with the following anisotropic Bernstein-type result (we refer to \cite{FC}, and to \cite{Dragos1} for more general anisotropic estimates):
\begin{lem}
\sl{There exists a constant $C>0$ such that for all function $f$, $\aa>0$, $1\leq q \leq p \leq \infty$ and all $0<r<R$, we have
\begin{equation}
\begin{cases}
\vspace{1mm}
\displaystyle{\|\chi (\frac{|D|}{R}) \chi (\frac{|D_3|}{r}) f\|_{L^p} \leq C \|f\|_{L^p},}\\
\displaystyle{\|\chi (\frac{|D|}{R}) \chi (\frac{|D_3|}{r}) f\|_{L^p} \leq C(R^2 r)^{\frac{1}{q}-\frac{1}{p}} \|\chi (\frac{|D|}{R}) \chi (\frac{|D_3|}{r}) f\|_{L^q}.}
\end{cases}
\end{equation}
Moreover if $f$ has its frequencies located in $\cC_{r,R}$, then
$$
\||D|^\aa f\|_{L^p} \leq C R^\aa \|f\|_{L^p}.
$$
}
\label{lemaniso}
\end{lem}

\subsection{Eigenvalues, projectors}

We begin with the eigenvalues and eigenvectors of matrix $\mathbb{B}(\xi, \ee)$. The main result of this section is the following proposition. We will only state the results and skip details as the proof is an adaptation of Proposition $3.1$ from \cite{FC5} (there in the anisotropic case).

\begin{prop}
\label{estimvp}
\sl{
There exists $\ee_0>0$ such that for all $\ee< \ee_0$, for all $r_\ee=\ee^m$ and $R_\ee =\ee^{-M}$, with $M<1/4$ and $3M+m<1$, and for all $\xi \in \mathcal{C}_{r_{\ee}, R_{\ee}}$, the matrix $\mathbb{B}(\xi, \ee) = \widehat{L-\frac{1}{\ee} \mathbb{P} \mathcal{A}}$ is diagonalizable and its eigenvalues have the following asymptotic expansions with respect to $\ee$:
\begin{equation}
\label{vp}
\begin{cases}
\vspace{0.2cm} \mu_0 = -\nu |\xi|^2,\\
\vspace{0.2cm} \mu = -(\nu\xi_1^2+ \nu\xi_2^2 + \nu'F^2\xi_3^2)\frac{|\xi|^2}{|\xi|_F^2} + \ee^2 D(\xi, \ee),\\
\vspace{0.2cm} \lambda = -\tau (\xi)|\xi|^2+i\frac{|\xi|_F}{\ee F   |\xi|}+ \ee E(\xi, \ee),\\
\vspace{0.2cm} \overline{\lambda} = -\tau (\xi)|\xi|^2-i\frac{|\xi|_F}{\ee F |\xi|}+ \ee \bar{E}(\xi, \ee),
\end{cases}
\end{equation}
where $|\xi|_F^2 = \xi_1^2 + \xi_2^2 + F^2\xi_3^3$, and $D,E$ denote remainder terms satisfying for all $\xi \in\cC_{r_\ee, R_\ee}$:
$$
\begin{cases}
\vspace{1mm}
\ee^2 |D(\xi, \ee)|\leq C_F |\nu-\nu'|^3 \ee^2 |\xi|^6\leq C_F |\nu-\nu'|^3\ee^{2-6M} \ll 1,\\
\ee |E(\xi, \ee)|\leq C_F |\nu-\nu'|^2 \ee |\xi|^4 \leq C_F |\nu-\nu'|^2\ee^{1-4M} \ll 1,\\
\ee |\n_\xi E(\xi, \ee)|\leq C_F |\nu-\nu'|^2 \ee |\xi|^3 \leq C_F |\nu-\nu'|^2\ee^{1-3M} \ll 1,\\
\end{cases}
$$
and
$$
\tau(\xi)=\frac{\nu}{2}\Big(1+\frac{F^2 \xi_3^2}{|\xi|_F^2}\Big) +\frac{\nu'}{2}\Big(1-\frac{F^2   \xi_3^2}{|\xi|_F^2}\Big) \geq \min (\nu,\nu')>0.
$$
Moreover, if we denote by $\mathcal{P}_i(\xi, \ee)$, the projectors onto the eigenspaces corresponding to $\mu$, $\lambda$ and $\overline{\lambda}$ ($i\in\{2,3,4\}$), and set
\begin{equation}
\label{projsevp} \mathbb{P}_i(u)=\mathcal{F}^{-1}\bigg(\mathcal{P}_i\big(\xi, \ee)(\widehat{u}(\xi)\big)\bigg),
\end{equation}
then for any divergence-free vector field $f$ whose Fourier transform is supported in $\mathcal{C}_{r_{\ee}, R_{\ee}}$, we have the following estimates:
\begin{equation}
\label{estimvp2}
\|\mathbb{P}_2 f\|_{H^s} \leq
\begin{cases}
{C_F \|f\|_{H^s}} & {\mbox{if} \quad \Omega(f)\neq 0,}\\
{C_F |\nu-\nu'| \ee^{1-(3M+m)} \|f\|_{H^s}} & {\mbox{if} \quad \Omega(f)= 0,}
\end{cases}
\end{equation}
and for $i=3,4$,
\begin{equation}
\label{estimvp34} \|\mathbb{P}_i f\|_{H^s}\leq C_F \ee^{-(m+M)} \|f\|_{H^s}.
\end{equation}
}
\end{prop}
\begin{rem}
\sl{We emphasize that the leading part of $\mu$ is the Fourier symbol of the quasi-geostrophic operator $\G$. Moreover, as we will see in what follows, the dispersion is related to the term $i\frac{|\xi|_F}{\ee F |\xi|}$, and when $F=1$ this term reduces to the constant $\frac{i}{\ee}$. This is why dispersion does not occur in the case $F=1$.}
\end{rem}

\subsection{Dispersion, Strichartz estimates}

The following result provides the Strichartz estimates satisfied by some projections of the solution of System\eqref{systdisp}:
\begin{prop}
\sl{Assume that $f$ satisfies \eqref{systdisp} on $[0,T[$ where $\div f_0=0$ and the frequencies of $f_0$ and $F$ are localized in $\cC_{r_\ee, R_\ee}$. Then there exists a constant $C_F>0$ such that for $i\in\{3,4\}$ and $p\geq 4$, we have
$$
\|\mathbb{P}_i f\|_{L_T^p L^\infty} \leq C_F \ee^{\frac{1}{p}-\left((\frac{5}{2}+\frac{4}{p})M +(2+\frac{4}{p})m)\right)} \left(\|f_0\|_{L^2} +\int_0^T \|\Fe(\tau)\|_{L^2} d\tau\right).
$$
\label{Strichartz1}
}
\end{prop}
\textbf{Proof: } We will only give a sketch of the proof and refer to \cite{FC,FC3,FC5} for details. First in the homogeneous case ($F=0$), the classical $TT*$ method allows us to write that:
$$
\|\mathbb{P}_i f\|_{L_T^p L^\infty} =\sup_{\psi \in \cB} \int_0^T \int_{\R^3} \mathbb{P}_i f(t,x) \psi(t,x) dx dt,
$$
where $\cB=\{\psi \in \cD(\R_+ \times \R^3), \quad \|\psi\|_{L_t^{\bar{p}} L^1} \leq 1\}$. Next, as $f$ has its frequencies localized in $\cC_{r_\ee, R_\ee}$ we obtain that:
\begin{multline}
\|\mathbb{P}_i f\|_{L_T^p L^\infty} =\sup_{\psi \in \cB} \int_0^T \int_{\R^3} \chi (\frac{|\xi|}{2R_\ee})(1-\chi (\frac{2|\xi_3|}{r_\ee})) e^{-t\tau (\xi)|\xi|^2+it\frac{|\xi|_F}{\ee F|\xi|}+ \ee tE(\xi, \ee)} \hat{\mathbb{P}_i f_0}(\xi) \hat{\psi}(t,\xi) d\xi dt.\\
\leq \|\mathbb{P}_i f_0\|_{L^2} \sup_{\psi \in \cB} \bigg(\int_0^T \int_0^T \|\psi(t)\|_{L^1} \|\psi(s)\|_{L^1} \|K(s,t,\ee, .)\|_{L^\infty} \bigg)^\frac{1}{2},
\label{estimTT}
\end{multline}
with
$$
K(s,t,\ee, x) =\int_{\R^3} e^{ix\cdot \xi} e^{-(t+s)\tau (\xi)|\xi|^2+i(t-s)\frac{|\xi|_F}{\ee F|\xi|}+ \ee tE(\xi, \ee) + \ee \bar{E}(\xi, \ee)} \chi (\frac{|\xi|}{2R_\ee})^2(1-\chi (\frac{2|\xi_3|}{r_\ee}))^2 dx.
$$
Adapting Section $4$ from \cite{FC5}, we obtain that:
$$
\|K(s,t,\ee, .)\|_{L^\infty} \leq C_F \frac{R_\ee^2}{r_\ee^2} \min\big(R_\ee, \frac{R_\ee^4}{r_\ee^2}\left(\frac{\ee}{|t-s|}\right)^{\frac{1}{2}}\big).
$$
\begin{rem}
\sl{We emphasize that the results from \cite{FC5} (Lemma $4.2$) are given for anisotropic viscosities and in the anisotropic space $L_{x_h,x_3}^{\infty,2}$, this is why the powers of $r_\ee,R_\ee$ are different in the present case. We also dropped the exponentials as we are in the finite time-case.}
\end{rem}
We then deduce that for all $\theta \in[0,1]$,
$$
\|K(s,t,\ee, .)\|_{L^\infty} \leq C_F \frac{R_\ee^{3(1+\theta)}}{r_\ee^{2(1+\theta)}} \left(\frac{\ee}{|t-s|}\right)^{\frac{\theta}{2}}.
$$
Using this for $\theta=4/p$ in \eqref{estimTT} (this is the reason why we need $p\geq 4$), we obtain the conclusion thanks to the Hardy-Littlewood theorem and the choice $r_\ee=\ee^m$ and $R_\ee=\ee^{-M}$ and using the projector estimate \eqref{estimvp34}. The non-homogeneous case easily follows. $\blacksquare$

From the previous property for $p=4$ we deduce the following result:
\begin{prop}
\sl{Under the same assumptions, there exists a constant $C_F>0$ such that for $i\in\{3,4\}$ and any $s\geq 0$, we have
$$
\|\mathbb{P}_i f\|_{L_T^4 B_{\infty, \infty}^s} +\|\mathbb{P}_i |D|^s f\|_{L_T^4 L^\infty} \leq C_F \ee^{\frac{1}{4}-\big(M(\frac{7}{2}+s)+3m\big)} \left(\|f_0\|_{L^2} +\int_0^T \|\Fe(\tau)\|_{L^2} d\tau\right).
$$
}
\label{Strichartz2}
\end{prop}
\textbf{Proof :} Thanks to the assumptions on the frequencies of $f_0$ and $F$, $f$ also has its frequencies in $\cC_{r_\ee, R_\ee}$ so that using Lemma \ref{lemaniso}, we get
$$
\|\mathbb{P}_i |D|^s f\|_{L_T^4 L^\infty} \leq C R_\ee^s \|\mathbb{P}_i f\|_{L_T^4 L^\infty},
$$
which immediately gives the second result thanks to Proposition \ref{Strichartz1}. For any $q\geq -1$ applying the Strichartz estimates to $\D_q f$ leads to:
\begin{multline}
\|\mathbb{P}_i \D_q f\|_{L_T^4 L^\infty} \leq C_F \ee^{\frac{1}{4}-\big(M(\frac{7}{2}+s)+3m\big)} \left(\|\D_q f_0\|_{L^2} +\int_0^T \|\D_q \Fe(\tau)\|_{L^2} d\tau\right)\\
\leq C_F \ee^{\frac{1}{4}-\big(M(\frac{7}{2}+s)+3m\big)} \left(\|f_0\|_{L^2} +\int_0^T \|\Fe(\tau)\|_{L^2} d\tau\right).
\end{multline}
Then for all $q\geq 0$, as $\D_q f$ is frequency localized on the ring $2^q \cC(0, \frac{3}{4}, \frac{8}{3})$ and in $\cC_{r_\ee, R_\ee}$, we have:
\begin{multline}
2^{qs} \|\mathbb{P}_i \D_q f\|_{L_T^4 L^\infty} \leq C \|\mathbb{P}_i |D|^s \D_q f\|_{L_T^8 L^\infty} \leq C R_\ee^s \|\mathbb{P}_i \D_q f\|_{L_T^8 L^\infty}\\
\leq C_F \ee^{\frac{1}{4}-\big(M(\frac{7}{2}+s)+3m\big)} \left(\|f_0\|_{L^2} +\int_0^T \|\Fe(\tau)\|_{L^2} d\tau\right).
\end{multline}
and for $q=-1$, we use that $2^{-s}\leq 1$. $\blacksquare$

\subsection{Final estimates}

The object of this section is to obtain the following estimates:
\begin{prop}
\sl{Assume that $\Ue$ solves \eqref{PE} on $[0,T]$ with the same divergence-free initial data as in Theorem \ref{thprincipal}. There exist $\ee_0>0$ and a positive constatnt $C$ such that for all $0<\ee <\ee_0$ and all $m,M>0$ with $M< 1/4$ and $3M+m <1$, and all $k\geq 0$ we have:
\begin{multline}
\||D|^k\Uosc\|_{L_T^4 L^\infty} \leq C_F \bigg(\ee^M +\ee^{\frac{m}{2}-M} +\ee^{1-(3M+m)} \bigg) T^\frac{1}{4}\|U_{0,\ee}\|_{H^{3+k}} e^{CV_\ee(T)}\\
+C_F \ee^{\frac{1}{4}-\big(M(5+k)+4m\big)} T^\frac{1}{2} \|U_{0,\ee}\|_{H^{3+k}}^2,
\end{multline}
where $V_\ee(t)=\int_0^t \|\n \Ue(\tau)\|_{L^\infty} d\tau$.
}
\label{estimUosc}
\end{prop}
\begin{rem}
\sl{In what follows, we will use this result in the case $k\leq 3$, $T\leq T_\ee^\gamma =\gamma \ln (\ln|\ln \ee|)$ and with $V_\ee(T) \leq -2\gamma \ln \ee$ and $\|U_{0,\ee}\|_{H^6} \leq C_0 \ee^{-5\beta}$. The previous estimates then turns into:
\begin{multline}
\||D|^k\Uosc\|_{L_T^4 L^\infty} \leq C_F \bigg(\ee^{M-(5\beta +2C\gamma)} +\ee^{\frac{m}{2}-(M+5\beta +2C\gamma)}\\ +\ee^{1-(3M+m+5\beta +2C\gamma)} +\ee^{\frac{1}{4}-\big(8M+4m+2(5\beta +2C\gamma)\big)}\bigg) (T_\ee^\gamma)^\frac{1}{2}.
\label{estimUoscTgamma}
\end{multline}
Moreover, under the following conditions:
\begin{equation}
\begin{cases}
5\beta +2C\gamma \leq \frac{M}{2},\\
M+(5\beta +2C\gamma) \leq \frac{m}{4},\\
3M+m+(5\beta +2C\gamma) \leq \frac{1}{2},\\
8M+4m+2(5\beta +2C\gamma) \leq \frac{1}{8},
\end{cases}
\end{equation}
which are satisfied for example if we have:
\begin{equation}
\begin{cases}
5\beta +2C\gamma \leq \frac{M}{2},\\
M\leq \frac{m}{4} \quad \mbox{and} \quad m\leq \frac{1}{44},
\end{cases}
\end{equation}
then we simply get, if $\ee$ is small enough,
\begin{equation}
\||D|^k\Uosc\|_{L_T^4 L^\infty} \leq C_F \ee^{\frac{M}{2}} (T_\ee^\gamma)^\frac{1}{2} \leq C_F \ee^{\frac{M}{4}}.
\label{estimUoscsimple}
\end{equation}
Note that these conditions on $m,M$ imply those in Proposition \ref{estimUosc} (which are required to use Proposition \ref{estimvp}).
}
\end{rem}
\textbf{Proof: } Let us cut the oscillating part into four parts as in \cite{FC}:
\begin{multline}
|D|^k\Uosc=\big(1-\chi (\frac{|D|}{R_\ee})\big)|D|^k\Uosc +\chi (\frac{|D|}{R_\ee})\chi (\frac{|D_3|}{r_\ee}) |D|^k\Uosc\\
+\cP_{r_\ee, R_\ee} |D|^k\mathbb{P}_2\Uosc +\cP_{r_\ee, R_\ee} |D|^k\mathbb{P}_{3+4}\Uosc\overset{def}{=} I_\ee+II_\ee+III_\ee+IV_\ee. 
\label{troncosc}
\end{multline}
Thanks to the injection $H^s(\R^3) \hookrightarrow L^\infty(\R^3)$ ($s>3/2$) and the apriori estimates from \eqref{estapriori} we roughly estimate the first term:
\begin{equation}
\|I_\ee\|_{L_T^4 L^\infty} \leq C T^\frac{1}{4} \|I_\ee\|_{L_T^\infty H^2} \leq \frac{C}{R_\ee} T^\frac{1}{4} \|\Uosc\|_{L_T^\infty H^{3+k}} \leq C T^\frac{1}{4}\ee^M \|\Ue\|_{L_T^\infty H^{3+k}}.
\label{estosc1}
\end{equation}
The second term is estimated thanks to Lemma \ref{lemaniso} (with $p=\infty$, $q=2$) and the Leray estimates:
\begin{equation}
\|II_\ee\|_{L_T^4 L^\infty} \leq C T^\frac{1}{4} (R_\ee^2 r_\ee)^\frac{1}{2} \||D|^k \Uosc\|_{L_T^\infty L^2} \leq C T^\frac{1}{4} \ee^{\frac{m}{2}-M} \|\Ue\|_{L_T^\infty H^k}.
\label{estosc2}
\end{equation}
As $\Om (\Uosc)=0$, we can use \eqref{estimvp2} to get:
\begin{multline}
\|III_\ee\|_{L_T^4 L^\infty} \leq C T^\frac{1}{4}\|III_\ee\|_{L_T^\infty H^2} \leq C_F T^\frac{1}{4}|\nu-\nu'| \ee^{1-(3M+m)} \|\cP_{r_\ee, R_\ee} |D|^k \Uosc\|_{L_T^\infty H^2}\\
\leq C_F T^\frac{1}{4}|\nu-\nu'| \ee^{1-(3M+m)} \|\Ue\|_{L_T^\infty H^{2+k}}.
\label{estosc3}
\end{multline}
The last term is estimated thanks to the Strichartz estimates from Proposition \ref{Strichartz2}:
$$
\|IV_\ee\|_{L_T^4 L^\infty} \leq C_F \ee^{\frac{1}{4}-\big(M(\frac{7}{2}+k)+3m\big)} \left(\|\cP_{r_\ee, R_\ee} U_{0,\ee,osc}\|_{L^2} +\int_0^T \|\cP_{r_\ee, R_\ee} \Fe(\tau)\|_{L^2} d\tau\right),
$$
where the external force term $\Fe$ is the right-hand side of \eqref{systosc}. As in \cite{FC} we roughly estimate each term from $\Fe$ (we refer to \eqref{systosc} for the decomposition) and obtain that (thanks to the Leray estimates and Lemma \ref{lemaniso}):
$$
\begin{cases}
\vspace{1mm}
\displaystyle{\|\cP_{r_\ee, R_\ee} F_1\|_{L_T^1 L^2} +\|\cP_{r_\ee, R_\ee} F_2\|_{L_T^1 L^2} \leq C R_\ee^\frac{3}{2} \int_0^T \|v_\ee\cdot \n \Ue\|_{L^1} d\tau \leq C R_\ee^\frac{3}{2} T^\frac{1}{2}\no^{-\frac{1}{2}} \|U_{0,\ee}\|_{L^2}^2,}\\
\vspace{1mm}
\displaystyle{\|\cP_{r_\ee, R_\ee} F_3\|_{L_T^1 L^2} \leq C_F r_\ee^{-1} R_\ee^\frac{3}{2} \int_0^T \|\n \Uosc \cdot \n \Ue\|_{L^1} d\tau \leq C_F r_\ee^{-1} R_\ee^\frac{3}{2} \no^{-1} \|U_{0,\ee}\|_{L^2}^2,}\\
\displaystyle{\|\cP_{r_\ee, R_\ee} F_4\|_{L_T^1 L^2} \leq  C_F |\nu-\nu'| R_\ee T^\frac{1}{2}\no^{-\frac{1}{2}} \|U_{0,\ee}\|_{L^2}.}
\end{cases}
$$
As $\|\cP_{r_\ee, R_\ee} U_{0,\ee,osc}\|_{L^2}\leq C r_\ee^{-1}\|U_{0,\ee,osc}\|_{\dot{H}^1}$ we finally obtain that the last term satisfies (to simplify we did not trace the viscosities and roughly estimated each norm by $\|U_{0,\ee}\|_{H^{3+k}}$.):
\begin{equation}
\|IV_\ee\|_{L_T^4 L^\infty} \leq C_F \ee^{\frac{1}{4}-\big(M(5+k)+4m\big)} T^\frac{1}{2} \|U_{0,\ee}\|_{H^{3+k}}^2.
\label{estosc4}
\end{equation}
Gathering \eqref{estosc1} to \eqref{estosc4} we end up with the desired result. $\blacksquare$

\section{$L^2$-estimates for the potential vorticity $\Ome$}

As emphasized in the second section when explaining the bootstrap argument, we need to obtain estimates for the potential vorticity in $L^2, L^\infty$ and $C^s(X)$. At first sight, the most natural way to do it seems, as in \cite{TH1, FC3}, to use transport-diffusion estimates on system \eqref{systomega}. We will be able to do this for the $L^\infty$-estimates, adaptating Proposition 1 from \cite{FCestimLp} which deals with the non-local operator $\G$.

Unfortunately, due to the first term from the right-hand side $(\nu-\nu') F\D \d_3 \theta_{\ee,osc}$, this method is useless for $p=2$ as we can only estimate this term by $|\nu-\nu'| F \|\n^3 \Ue\|_{L^2}$ that we do not control (the best we can hope for is to estimate this by a negative power of $\ee$). In order to overcome this problem we need to bound $\|\Uqg\|_{\dot{H}^1}$. Thanks to Proposition \ref{propdecomposcqg} and the fact that $\Uosc$ is small, we will simply study $\|\Ue\|_{\dot{H}^1}$ and then use that $\|\Ome\|_{L^2} \leq C\|\Uqg\|_{\dot{H}^1} \leq C\|\Ue\|_{\dot{H}^1}$.

Another motivation to estimate directly $\|\Ue\|_{\dot{H}^1}$ is that getting $L^p$-estimates (see Proposition \ref{appendicestimLp} in the appendix) requires a $L_t^\infty L^6$-bound for the transport term $v_{\ee}$ that will be a simple consequence of the result on $L_t^\infty \dot{H}^1$ thanks to the Sobolev injection.

\subsection{$\dot{H}^1$-estimates for the solution $\Ue$}

The aim of this section is to prove:
\begin{prop}
\sl{Under the assumptions of Theorem \ref{thprincipal}, and with the notations from \eqref{bootstrap}, there exists $C_F,M>0$ and $\ee_0$ such that if $5\beta +2C\gamma \leq M/6$, then for all $0<\ee \leq \ee_0$ and all $t\leq \min(T_\ee, T_\ee^\gamma)$, we have:
\begin{equation}
\frac{1}{2} \|\Ue(t)\|_{\dot{H}^1}^2 +\no \int_0^t \|\n \Ue(\tau)\|_{\dot{H}^1}^2 d\tau \leq C_F
\label{estimUH}
\end{equation}
}
\label{PropUH}
\end{prop}
\textbf{Proof :} we emphasize that the quasi-geostrophic/oscillating decomposition is an orthogonal decomposition of the solution so that we have:
$$
\|\Ue\|_{\dot{H}^1}^2 \sim \|\Uqg\|_{\dot{H}^1}^2 +\|\Uosc\|_{\dot{H}^1}^2 \sim \|\Ome\|_{L^2}^2 +\|\Uosc\|_{\dot{H}^1}^2 ,
$$
and to obtain the desired estimates we could, as in \cite{Chemin2} estimate separatedly the potential vorticity and the oscillating part. But in our case, as $F\neq 1$, it will be easier to use dispersion phenomena together with the quasi-geostrophic structure and estimate the following (equivalent) norm:
$$
\|\Ue\|_{\dot{H}_F^1}^2 \overset{def}{=} -(\D_F \Ue|\Ue)_{L^2}.
$$
Taking the corresponding $\dot{H}_F^1$-inner product of \eqref{PE} with $\Ue$, we obtain that
$$
\frac{1}{2} \frac{d}{dt} \|\Ue\|_{\dot{H}_F^1}^2 -(L \Ue| \Ue)_{\dot{H}_F^1}=-(v_\ee\cdot \n\Ue| \Ue)_{\dot{H}_F^1}.
$$
As usual $-(L \Ue| \Ue)_{\dot{H}_F^1} \geq \no \|\n \Ue\|_{\dot{H}_F^1}^2$ and thanks to the fact that $\Ue=\Uosc+\Uqg$ we develop the right-hand side as follows:
\begin{multline}
(v_\ee\cdot \n\Ue| \Ue)_{\dot{H}_F^1} =(v_\ee\cdot \n\Ue| \Uosc)_{\dot{H}_F^1} +(v_\ee\cdot \n\Uosc| \Uqg)_{\dot{H}_F^1} +(v_{\ee,osc}\cdot \n\Uqg| \Uqg)_{\dot{H}_F^1}\\
+(v_{\ee,QG}\cdot \n\Uqg| \Uqg)_{\dot{H}_F^1} =B_1 +B_2+ B_3 +B_4.
\end{multline}
What is more unusual is that the most dangerous term, namely $B_4$ (i.-e. the only one that may be large as it does not involve the oscillating part, we will show that the three other terms are small), is equal to zero. We emphasize that if we had computed the classical $\dot{H}^1$-innerproduct it was not true anymore (as for any Navier-Stokes type system), and this term would have obstructed any use of the logarithmic (or Gronwall) estimates. To show this we simply use the following elementary computation related to the quasi-geostrophic decomposition: for any function $f$, we have:
$$
(f|\Uqg)_{\dot{H}_F^1} =-(f|\left(
\begin{array}{c}
-\d_2\\
\d_1\\
0\\
-F\d_3
\end{array}
\right)\Ome)_{L^2} =(\Om(f)|\Ome)_{L^2}.
$$
Then, thanks to this, point $6$ from Proposition \ref{propdecomposcqg} and the fact that $\div v_{\ee,QG}=0$, we obtain:
$$
B_4=(v_{\ee,QG}\cdot \n\Uqg| \Uqg)_{\dot{H}_F^1}=(\Om(v_{\ee,QG}\cdot \n\Uqg)| \Ome)_{L^2} =(v_{\ee,QG}\cdot \n\Ome| \Ome)_{L^2}=0.
$$
Writing
$$
\begin{cases}
|B_1| \leq C_F \|v_\ee\|_{L^2} \|\n \Ue\|_{L^2} \|\n^2 \Uosc\|_{L^\infty},\\
|B_2| \leq C_F \|v_\ee\|_{L^2} \|\n^2 \Ue\|_{L^2} \|\n \Uosc\|_{L^\infty},\\
|B_3| \leq C_F \|v_{\ee, osc}\|_{L^\infty} \|\n \Uqg\|_{L^2} \|\n^2 \Uqg\|_{L^2},\\
\end{cases}
$$
thanks to the Leray estimates \eqref{estimLeray}, the a priori estimates \eqref{estapriori} and the fact that $\cP,\cQ$ are pseudo-differential operators of order zero, we get that for all $t\leq \min (T_\ee^*, T_\ee)$ :
\begin{multline}
\frac{1}{2} \|\Ue(t)\|_{\dot{H}_F^1}^2 +\no \int_0^t \|\n \Ue(\tau)\|_{\dot{H}_F^1}^2 d\tau\\
\leq \frac{1}{2} \|U_{0,\ee}\|_{\dot{H}_F^1}^2 +C_F \int_0^t \|\Ue(\tau)\|_{H^6}^2 \bigg(\|\Uosc\|_{L^\infty} +\|\n \Uosc\|_{L^\infty} +\|\n^2 \Uosc\|_{L^\infty} \bigg) d\tau\\
\leq C_F \bigg(\|U_{0,\ee}\|_{\dot{H}^1}^2 +\|U_{\ee, 0}\|_{H^6}^2 e^{2CV_\ee(t)} t^\frac{3}{4} \big(\|\Uosc\|_{L_t^4 L^\infty} +\|\n \Uosc\|_{L_t^4 L^\infty} +\|\n^2 \Uosc\|_{L_t^4 L^\infty}\big) \bigg).
\end{multline}
Therefore, thanks to Proposition \ref{estimUosc} and \eqref{estimUoscTgamma} (we roughly use it even if we do not have derivatives of order three here) and as $t\leq \min(T_\ee,\gamma \ln (\ln|\ln \ee|))$, $V_\ee(t) \leq K_\ee= -2\gamma \ln \ee$ and $\|U_{0,\ee}\|_{H^6} \leq C_0 \ee^{-5\beta}$ we have:
\begin{multline}
\frac{1}{2} \|\Ue(t)\|_{\dot{H}_F^1}^2 +\no \int_0^t \|\n \Ue(\tau)\|_{\dot{H}_F^1}^2 d\tau\\
\leq C_F \bigg(\|U_{0,\ee}\|_{\dot{H}^1}^2 +(T_\ee^\gamma)^\frac{5}{4} \big( \ee^{M-3(5\beta +2C\gamma)} +\ee^{\frac{m}{2}-\big(M+3(5\beta +2C\gamma)\big)}\\
+\ee^{1-\big(3M+m+3(5\beta +2C\gamma)\big)} +\ee^{\frac{1}{4}-\big(8M+4m+4(5\beta +2C\gamma)\big)}\big) \bigg).
\end{multline}
So, similarly, under the following conditions:
\begin{equation}
\begin{cases}
3(5\beta +2C\gamma) \leq \frac{M}{2},\\
M+3(5\beta +2C\gamma) \leq \frac{m}{4},\\
3M+m+3(5\beta +2C\gamma) \leq \frac{1}{2},\\
8M+4m+4(5\beta +2C\gamma) \leq \frac{1}{8},
\end{cases}
\end{equation}
which are implied for example if we have:
\begin{equation}
\begin{cases}
5\beta +2C\gamma \leq \frac{M}{6},\\
M\leq \frac{m}{6} \quad \mbox{and} \quad m\leq \frac{1}{49},
\end{cases}
\label{condmM}
\end{equation}
then as $\min (\frac{M}{2}, \frac{m}{4}, \frac{1}{2}, \frac{1}{8}) =\frac{M}{2}$, if $\ee$ is small enough we can write that:
\begin{multline}
\frac{1}{2} \|\Ue(t)\|_{\dot{H}_F^1}^2 +\no \int_0^t \|\n \Ue(\tau)\|_{\dot{H}_F^1}^2 d\tau \leq C_F \big(\|U_{0,\ee}\|_{\dot{H}^1}^2 +(T_\ee^\gamma)^\frac{5}{4} \ee^{\frac{M}{2}} \big)\\
\leq C_F \big( C_0 +\ee^{\frac{M}{4}}) \leq C_F.
\end{multline}
This concludes the proof. $\blacksquare$

\subsection{$L^2$-estimate for the potential vorticity}

As a direct consequence of \eqref{estimUH} we obtain the following result:
\begin{prop}
\sl{Under the assumptions of Theorem \ref{thprincipal}, and with the notations from \eqref{bootstrap}, there exists $C_F,M>0$ (with $M<1/294$) and $\ee_0$ such that if $5\beta +2C\gamma \leq M/6$, then for all $0<\ee \leq \ee_0$ and all $t\leq \min(T_\ee, T_\ee^\gamma)$, we have:
\begin{equation}
\frac{1}{2} \|\Ome(t)\|_{L^2}^2 +\no \int_0^t \|\n \Ome(\tau)\|_{L^2}^2 d\tau \leq C_F
\label{estimOmegaL2}
\end{equation}
}
\end{prop}

\section{Proof of the main result}

\subsection{$L^\infty$-estimates for the potential vorticity}

The object of this section is to obtain $L^\infty$-estimates for the potential vorticity. This result is very close to Proposition $1$ from \cite{FCestimLp}.

\begin{prop}\sl{Under the assumptions of Proposition \ref{PropUH}, and with the notations from \eqref{bootstrap}, there exists $C_F, D_F,  C_0',M>0$ and $\ee_0$ ($D_F$ depending on $F,\nu$ and $\nu'$, $C_0'$ depending on $\|\Om_0\|_{L^\infty}$) such that if $5\beta +2C\gamma \leq M/6$, then for all $0<\ee \leq \ee_0$ and all $t\leq \min(T_\ee, T_\ee^\gamma)$, we have:
\begin{equation}
\begin{cases}
\vspace{1mm}
\displaystyle{\frac{1}{2} \|\Ue(t)\|_{\dot{H}^1}^2 +\no \int_0^t \|\n \Ue(\tau)\|_{\dot{H}^1}^2 d\tau \leq C_F,}\\
\|\Ome\|_{L_t^\infty L^\infty}\leq C_0' e^{t D_F}.
\end{cases}
\label{estimOmLinf}
\end{equation}
}
\end{prop}
\textbf{Proof :} the $\dot{H}^1$-estimates were obtained in the previous section and the ideas to get the $L^\infty$-estimates are the same as in \cite{FCestimLp} so we will skip details and only focus on what is different. Let us define for all $\xi\neq 0$, $q_0(\xi)$ by
$$
\hat{\G u}(\xi)= -\frac{|\xi|^2}{|\xi|_F^2} (\nu \xi_1^2+ \nu \xi_2^2 +\nu' F^2 \xi_3^2) \hat{u}(\xi) \overset{\mbox{def}}{=}-q_0(\sqrt{\no}\xi)\hat{u}(\xi).
$$
If we denote $M=\frac{\nu}{\no}$ and $M'=\frac{\nu'}{\no}$, (recall that $\nu_0=\min(\nu, \nu')>0$) then
$$
\begin{cases}
\min(M, M')=1,\\
\max(M,M')=M_{visc}=\frac{\max(\nu, \nu')}{\min(\nu, \nu')},
\end{cases}
$$
so that
$$
q_0(\xi)= \frac{|\xi|^2}{|\xi|_F^2} (M \xi_1^2+ M \xi_2^2 +M' F^2 \xi_3^2) \geq |\xi|^2.
$$
And, as explained in \cite{FCestimLp}, we can write that $e^{t\G} u= K_t * u$, where the kernel $K_t(x)$ is defined for all $t,x$ by:
$$
K_t(x)=\frac{1}{\sqrt{\no t}^3} K_1(\frac{x}{\sqrt{\no t}}), \quad \mbox{with} \quad K_1(x)=\frac{1}{(2\pi)^3} \int_{\R^3} e^{ix\cdot \xi} e^{-q_0(\xi)} d\xi.
$$
We recall that in \cite{FCestimLp} we obtained there exists a constant $C_F'>0$ depending on $F,M_{visc}$ such that
\begin{equation}
\|K_1\|_{L^1} +\|\n K_1\|_{L^\frac{6}{5}} \leq C_F'.
\label{estimnoyau}
\end{equation}
Thanks to the Duhamel form, we obtain that for all $t\leq \min(T_\ee, T_\ee^\gamma)$,
$$
\Ome(t)= e^{t\G} \Om_{0,\ee} +\int_0^t e^{(t-\tau)\G} \bigg(-\div(v_\ee\otimes \Ome)(\tau)+(\nu-\nu')F\d_3 \D \theta_{\ee,osc}(\tau) +q_\ee(\tau)\bigg) d\tau.
$$
Then as in the cited paper, replacing $q_\ee$ by its value and thanks to \eqref{estimnoyau} and convolution estimates,
\begin{multline}
\|\Ome(t)\|_{L^\infty}\\
 \leq C_F'\left( \|\Om_{0,\ee}\|_{L^\infty} +\int_0^t \big(F|\nu-\nu'|\|\d_3 \D \theta_{\ee,osc}(\tau)\|_{L^\infty} +\|\n\Ue(\tau)\cdot \n \Uosc(\tau)\|_{L^\infty} \big) d\tau\right)\\
+\int_0^t \frac{1}{\sqrt{\no(t-\tau)}^4} \sqrt{\no(t-\tau)}^{3\cdot \frac{5}{6}} \|\n K_1\|_{L^\frac{6}{5}} \|v_\ee(\tau)\|_{L^6} \|\Ome(\tau)\|_{L^\infty} d\tau.
\end{multline}
We recall that as assumed in Theorem \ref{thprincipal}, $\Ome =\ee^{-3\beta} h(\ee^{\beta}.) *\Om_0$. Then, using the injection $H^2(\R^3) \hookrightarrow L^\infty(\R^3)$ we get:
\begin{multline}
\|\Ome(t)\|_{L^\infty} \leq C_F'\left( \|\Om_0\|_{L^\infty} +(T_\ee^\gamma)^\frac{3}{4} \big(\|\n^3 \Uosc\|_{L_t^4 L^\infty} +\|U_{0,\ee}\|_{H^6} e^{CV_\ee(t)}\|\n \Uosc\|_{L_t^4 L^\infty} \big) \right)\\
+C_F' \no^{-\frac{3}{4}} \big(\int_0^t (t-\tau)^{-\frac{3}{4}}d\tau\big) \|v_\ee\|_{L^\infty \dot{H}^1} \|\Ome\|_{L_t^\infty L^\infty}.
\end{multline}
and thanks to \eqref{estimUoscTgamma} and \eqref{estimUH}, 
\begin{multline}
\|\Ome\|_{L_t^\infty L^\infty} \leq C_F'\bigg( \|\Om_0\|_{L^\infty} +(T_\ee^\gamma)^\frac{5}{4} \big( \ee^{M-2(5\beta +C\gamma)} +\ee^{\frac{m}{2}-2(M+5\beta +2C\gamma)}\\
+\ee^{1-\big(3M+m+2(5\beta +2C\gamma)\big)} +\ee^{\frac{1}{4}-\big(8M+4m+3(5\beta +2C\gamma)\big)}\big) \bigg) +4C_F' \no^{-\frac{3}{4}} t^\frac{1}{4}  C_F'\|\Ome\|_{L_t^\infty L^\infty}.
\end{multline}
The conditions on $M,m,\beta, \gamma$ under which the exponents of $\ee$ are positive are implied by \eqref{condmM}, then there exists $\ee_0>0$ (only depending on $\beta, \gamma$) such that if $\ee\leq \ee_0$ then
\begin{multline}
(T_\ee^\gamma)^\frac{5}{4} \bigg( \ee^{M-2(5\beta +C\gamma)} +\ee^{\frac{m}{2}-2(M+5\beta +2C\gamma)} +\ee^{1-\big(3M+m+2(5\beta +2C\gamma)\big)}\\
+\ee^{\frac{1}{4}-\big(8M+4m+3(5\beta +2C\gamma)\big)}\bigg) \leq 1
\end{multline}
and if $t$ is so small that $4C_F' \no^{-\frac{3}{4}} t^\frac{1}{4}  C_F' \leq \frac{1}{2}$, that is
\begin{equation}
t\leq \frac{\no^3}{(8C_F')^4},
\end{equation}
then we obtain that
$$
\|\Ome\|_{L_t^\infty L^\infty} \leq C_F'(\|\Om_0\|_{L^\infty} +1).
$$
As we want a result for large times $t$, we will globalize the estimates with the usual method: if we subdivide the interval $[0,t]$ into $0=T_0<T_1<...<T_N=t$ such that for any $i\in\{0,...,N-1\}$,
\begin{equation}
T_{i+1}-T_i\sim \frac{\no^3}{(8C_F')^4},
\label{condTi}
\end{equation}
then the previous arguments imply that if $0<\ee\leq \ee_0$ for any $i\in\{0,...,N-1\}$,
$$
\|\Ome\|_{L\infty([T_i,T_{i+1}], L^\infty)} \leq C_F'(\|\Ome(T_i)\|_{L^\infty} +1) \leq C_0' (C_F')^i.
$$
Next summing \eqref{condTi} for $i\in\{0,...,N-1\}$, we obtain that $N\sim \frac{(8C_F')^4}{\no^3} t$ so that finally if we denote $D_F \overset{def}{=} \frac{(8C_F')^4}{\no^3} \ln C_F'$ we obtain the desired estimate for $\Ome$. $\blacksquare$

\subsection{Tangential regularity}

\subsubsection{Advected family of vector fields}

We refer to the appendix for an introduction to the notations related to the vortex patches. As mentionned before, a crucial ingredient in the proof of the main result of this article is to use the logarithmic estimates \eqref{estimlog}. However we need to be careful that if $\Om_0$ is $C^s(X_0)$ where $X_0=\{X_{0,\lambda}, \lambda=1,...,N\}$ is a fixed admissible system of $C^s$-vectorfields, we will not measure $\|\Ome\|_{C^s(X_0)}$ but $\|\Ome\|_{C^s(\Xtl)}$ where $\Xtl$ is the solution of the following transport equation:
\begin{equation}
\begin{cases}
\d_t \Xtl+ v_\ee.\n \Xtl = \Xtl\cdot \n v_\ee\\
{\Xtl}_{/t=0}= X_{0, \lambda}.
\end{cases}
\label{systXtl}
\end{equation}
The regularity is preserved by this transformation : we refer to \cite{GSR} for the proof of the fact that $X_{0, \lambda} \in
C^s \Rightarrow \Xtl\in C^s$. We refer to \cite{Chemin1, GSR, TH1, Dutrifoy1} for more details about the persistence of the tangential (or stratified) regularity in the vortex patches theory and we recall that we denote for any function $w$ (using the same notation as in \cite{Dutrifoy1, FC3}):
$$
\Xtl(x,D)w =\div (w \otimes \Xtl).
$$
Then in our case, \eqref{estimlog} becomes
\begin{equation}
\|\n \Uqg\|_{L^{\infty}} \leq C\Big(\|\Ome\|_{L^2}+ \|\Ome\|_{L^{\infty}} \log \big(e+ \frac{\|\Ome\|_{C^s(\Xtl)}}{\|\Ome\|_{L^{\infty}}} \big) \Big),
\label{estimlog2}
\end{equation}
where
\begin{equation}
\|\Ome\|_{C^s(\Xtl)} =\|\Ome\|_{L^{\infty}}+\|[\Xtl]^{-1}\|_{L^{\infty}} +\sum_{\lambda=1}^N \big( \|\Xtl\|_{C^s}+ \|\Xtl(x, D)\Ome\|_{C^{s-1}} \big).
\end{equation}
and there remains for us to estimate the last three terms. Only the last one requires a careful study as the other terms are strictly the same as in \cite{FC3}. For now we will estimate these other terms and state in the next section the smoothing effect provided by a priori estimates on $\Ome$ that will help us estimating $\Xtl(x, D)\Ome$. In what follows we will first obtain estimates in a small time interval $[T_1,T_2]$ and then globalize the results at the end of the article.
\begin{rem}
\sl{For more simplicity, in this section we will denote $X^\ee(t)$ instead of $\Xtl$.}
\end{rem}
System \eqref{systXtl} has been extensively studied (see \cite{Chemin1, GSR, TH1, Dutrifoy1, FC3}) and we refer to \cite{FC3} for the following estimate: there exists a constant $C_s>0$ such that for all $t\in [T_1, T_2]$,
\begin{multline}
\|X^\ee(t)\|_{C^s} \leq \|X^\ee(T_1)\|_{C^s} +C_s \int_{T_1}^t \|X^\ee(t')\|_{C^s}  \|\Ue\|_{Lip}dt'\\
+C_s \int_{T_1}^t \|\Xdo(t')\|_{C^{s-1}} dt'+ C_s \int_{T_1}^t \|X^\ee(t')\|_{L^{\infty}} \|\Uosc\|_{C^{s+1}}dt'.
\end{multline}
Without any change, as in \cite{FC3} we refer to \cite{Chemin1, GSR} for the fact that for all $t\in[0,T_\ee[$,
\begin{equation}
\begin{cases}
\vspace{1mm}
\displaystyle{\|X^\ee(t)\|_{L^{\infty}} \leq \|X_0\|_{L^{\infty}} e^{CV_\ee(t)} \leq C_0 \ee^{-2C\gamma},}\\
\displaystyle{\|[\Xtl]^{-1}\|_{L^{\infty}} \leq C_0 e^{CV_\ee(t)} \leq C_0 \ee^{-2C\gamma}.}
\end{cases}
\label{estimXe}
\end{equation}
Using this and estimates \eqref{estimUoscTgamma} and \eqref{estimUoscsimple}, we obtain that under condition \eqref{condmM} on $m,M,\beta$ and $\gamma$ there exists $\ee_0>0$ such that if $0<\ee\leq \ee_0$ we have,
$$
C_s \int_{T_1}^t \|X^\ee(t')\|_{L^{\infty}} \|\Uosc\|_{C^{s+1}}dt' \leq C_{F,s} \ee^\frac{M}{2} (T_\ee^\gamma)^\frac{1}{2} \leq 1.
$$
and therefore:
\begin{multline}
\|X^\ee(t)\|_{C^s} \leq \|X^\ee(T_1)\|_{C^s} +1 +C_s \int_{T_1}^t \|X^\ee(t')\|_{C^s} \|\Ue\|_{Lip}dt'\\
+C_s \int_{T_1}^t \|\Xdo(t')\|_{C^{s-1}} dt'.
\label{estimXe2}
\end{multline}

\subsubsection{Smoothing effect for $\Ome$}

The aim of this section is to state the heat regularization occuring for the potential vorticity $\Ome$ in Besov and H\"older spaces. We recall that $\Ome$ satisfies system \eqref{systomega}. Thanks to Proposition \ref{PropUH}, under conditions \eqref{condmM} and if $0<\ee\leq \ee_0$, we have
\begin{equation}
\|\ve\|_{L^\infty([T_1,T_2],L^6)} \leq \|\Ue\|_{L^\infty([T_1,T_2],\dot{H}^1)}\leq C_F.
\label{estimvee}
\end{equation}
If in addition $T_1\leq T_2$ are such that:
\begin{equation}
\begin{cases}
2C_F(T_2-T_1)^{\frac{1}{4}} \leq \no^3,\\
e^{C\int_{T_1}^{T_2} \|\nabla \Ue (\tau)\|_{L^\infty}}-1 \leq \frac{1}{C_F M_{visc}},
\end{cases}
\label{condT1T2b}
\end{equation}
then from Theorem \ref{thC2} (see Appendix, this result is proved in \cite{FCestimLp}), for all $t\in[T_1,T_2]$
\begin{multline}
\|\Ome\|_{\tilde{L}^\infty([T_1,t], C_*^0)} +\no \|\Ome\|_{\tilde{L}^1([T_1,t], C_*^2)}\\
\leq C_F \left[\|\Ome(T_1)\|_{L^\infty} +\int_{T_1}^t \bigg(F|\nu-\nu'|\|\d_3 \D \theta_{\ee,osc}(\tau)\|_{L^\infty} +\|q_\ee(\tau)\|_{L^\infty}\bigg)d\tau \right].
\end{multline}
Thanks to \eqref{estimUoscTgamma} and \eqref{estimUoscsimple}, we estimate as usual the oscillating part and obtain:
$$
\|\Ome\|_{\tilde{L}^\infty([T_1,t], C_*^0)} +\no \|\Ome\|_{\tilde{L}^1([T_1,t], C_*^2)} \leq C_F \left(\|\Ome(T_1)\|_{L^\infty} +\ee^\frac{M}{2} (T_\ee^\gamma)^\frac{1}{2} (t-T_1)^\frac{3}{4}\right).
$$
The last term is small so there exists $\ee_0>0$ such that for all $0<\ee \leq \ee_0$ it is bounded by $1$, then using \eqref{estimOmLinf}, we end up with:
\begin{equation}
\|\Ome\|_{\tilde{L}^\infty([T_1,t], C_*^0)} +\no \|\Ome\|_{\tilde{L}^1([T_1,t], C_*^2)} \leq C_F (e^{T_1 D_F}+1) \leq 2C_Fe^{T_1 D_F}.
\label{estimOm2}
\end{equation}

\subsubsection{Study of $\Xtl(x, D)\Ome$}

As explained before, we now turn to the most difficult term. We begin with the transport-diffusion system satisfied by $\Xtl(x, D)\Ome$: similarly to \cite{FC3} (the differences are the non-local operator and the additionnal term $F(\nu-\nu')\d_3 \D \theta_\ee$) we easily obtain that:
\begin{multline}
(\d_t+ v_\ee\cdot \n -\G) \big(\Xtl(x, D) \Ome\big)\\
=(\nu-\nu')F \Xtl(x,D) \big(\d_3 \D \theta_\ee\big) +\Xtl(x,D) q_\ee+[\Xtl(x,D),\G] \Ome,\\
\label{transportregtg}
\end{multline}
with the initial data
$$
{\Xtl(x,D)\Ome}_{/t=0} =X_{0, \lambda}(x, D) \Om_{0, \ee}.
$$
As in \cite{TH1, FC3} all the difficulty will lie in estimating the last term $[\Xtl(x,D),\G] \Ome$ and in our case, as in \cite{FCestimLp}, everything will be far more tricky as we have to deal with the non-local diffusion operator $\G$ instead of the Laplacian. As was first observed in \cite{TH1} this commutator is a sum of terms of different regularity. From \eqref{decpgamma} that provides the decomposition of $\G$ into its local and purely non-local parts (we refer to \cite{FCestimLp} for more details) we can decompose the commutator into:
\begin{multline}
[\Xtl(x,D),\G]\Ome =\A +\B\\
=[\Xtl(x,D),\GL]\Ome +(\nu-\nu') F^2(1-F^2) [\Xtl(x,D), \La^2]\Ome.
\end{multline}
The first term will be decomposed exactly as in \cite{TH1, FC3} into:
\begin{equation}
\A =F_\A + G_\A,
\label{decomutA}
\end{equation}
with
$$
F_\A = -\sum_{i=1}^3 \bigg( \d_i R(\Ome, \GL \Xtli) +2\d_i R(^t\n \Xtli, M_{\nu, \nu', F} \n \Ome)\bigg),
$$
and
\begin{multline}
G_\A = -\sum_{i=1}^3 \bigg( \d_i T_{\Ome} \GL \Xtli +\d_i T_{\GL \Xtli} \Ome\\
+2\d_i T_{^t\n \Xtli} M_{\nu, \nu', F} \n \Ome +2\d_i T_{M_{\nu, \nu', F} \n \Ome} \n \Xtli\bigg),
\end{multline}
where $M_{\nu, \nu', F}$ is the diagonal matrix defined by
$$
\left(\begin{array}{ccc}
\nu & 0 & 0\\
0 & \nu & 0\\
0 & 0 & (1-F^2)\nu +F^2 \nu'
\end{array}\right),
$$
and $R$ and $T$ correspond to the Bony decomposition and are defined in \eqref{Bonydecp} (we refer to \cite{Bony, Chemin1, Dbook} for precise studies of these operators). The second term first needs to be rewritten into :
\begin{multline}
\B =-(\nu-\nu') F^2(1-F^2)\left(\Xtl(x,D) (\La^2\Ome) -\La^2 \big(\Xtl(x,D)\Ome\big)\right)\\
=-(\nu-\nu') F^2(1-F^2)\sum_{i=1}^3 \d_i\bigg(\Xtli \La^2\Ome -\La^2(\Xtli \Ome)\bigg)\\
=-(\nu-\nu') F^2(1-F^2)\sum_{i=1}^3 \d_i\bigg( \big(\Xtli \La\big(\La\Ome\big) -\La(\Xtli \La\Ome) \big) +\La\big(\Xtli \La\Ome-\La(\Xtli \Ome)\big)\bigg).
\end{multline}
And using the bilinear operator $M$ defined in \eqref{defM}, we can decompose $\B$ into:
\begin{multline}
\B =-(\nu-\nu') F^2(1-F^2)\sum_{i=1}^3 \d_i\bigg[ \La\Xtli \La\Ome +M(\Xtli, \La \Ome) +\La\bigg(\La\Xtli \Ome +M(\Xtli, \Ome)\bigg)\bigg]\\
=F_\B + G_\B,
\label{decomutB}
\end{multline}
with
\begin{multline}
F_\B =-(\nu-\nu') F^2(1-F^2)\sum_{i=1}^3 \d_i\bigg[ R\big(\La\Xtli, \La\Ome\big) +\M\bigg(R\big(\Xtli ,\La \Ome\big)\bigg)\\
+\La R\big(\La\Xtli, \Ome\big) + \La \M\bigg(R\big(\Xtli, \Ome\big)\bigg)\bigg],
\end{multline}
and
\begin{multline}
G_\B=-(\nu-\nu') F^2(1-F^2)\sum_{i=1}^3 \d_i\bigg[ T_{\La\Xtli} \La\Ome +T_{\La \Ome} \La\Xtli +\M(T_{\Xtli} \La \Ome) +\M(T_{\La \Ome} \Xtli)\\
+\La\bigg(T_{\La\Xtli} \Ome +T_{\Ome} \La\Xtli +\M\big(T_{\Xtli} \Ome\big) +\M\big(T_{\Ome} \Xtli\big) \bigg)\bigg],
\end{multline}
where we denote for any smooth functions $f,g$:
$$
\begin{cases}
\displaystyle{\M(T_f g) =\sum_{q\geq -1} M(S_{q-1} f, \D_q g),}\\
\displaystyle{\M\big(R(f,g)\big) =\sum_{q\geq -1} \sum_{\aa =-1}^1 M(\D_q f, \D_{q+\aa} g).}
\end{cases}
$$
As in \cite{TH1, FC3} we decompose the commutators into the sum of two terms which have different regularity. This is the object of the following result:
\begin{prop}
\sl{If $s\in]0,1[$ there exists a constant $C$ such that for all $0\leq T_1 \leq T_2\leq \min(T_\ee, T_\ee^\gamma)$ satisfying condition \eqref{condT1T2b}, we have for all $t\in[T_1,T_2]$:
\begin{equation}
\begin{cases}
\vspace{1mm}
\displaystyle{\|F_\A\|_{\tilde{L}^1([T_1,t], C^{s-1})} +\|F_\B\|_{\tilde{L}^1([T_1,t], C^{s-1})} \leq C_{F,s} \max(\nu, \nu') \|\Xtol\|_{L^\infty([T_1,t],C^s)} \|\Ome\|_{\tilde{L}^1([T_1,t], C_*^2)},}\\
\displaystyle{\|G_\A(t)\|_{C^{s-3}} +\|G_\B(t)\|_{C^{s-3}} \leq C_{F,s} \max(\nu, \nu') \|\Ome(t)\|_{L^\infty} \|\Xtl\|_{C^s}.}
\end{cases}
\end{equation}
}
\end{prop}
This immediately implies, thanks to \eqref{estimOmLinf} and \eqref{estimOm2} that for all $t\in[T_1,T_2]$:
\begin{equation}
\begin{cases}
\vspace{1mm}
\displaystyle{\|F_\A\|_{\tilde{L}^1([T_1,t], C^{s-1})} +\|F_\B\|_{\tilde{L}^1([T_1,t], C^{s-1})} \leq C_F M_{visc} e^{T_1 D_F} \|\Xtol\|_{L^\infty([T_1,t],C^s)},}\\
\displaystyle{\frac{1}{\no}\bigg(\|G_\A(t)\|_{C^{s-3}} +\|G_\B(t)\|_{C^{s-3}}\bigg) \leq C_F M_{visc} e^{T_1 D_F}\|\Xtl\|_{C^s}.}
\end{cases}
\label{estimFG}
\end{equation}
\textbf{Proof :} even if the estimates of $F_\A$ and $G_\A$ are obtained exactly as in \cite{TH1, FC3} we will give a quick proof of them as it helps understanding how to deal with $F_\B$ and $G_\B$ which are more tricky. Let us begin with $F_\A$. In this proof, in order to simplify, we will drop the summation in $i$:
\begin{multline}
\|F_\A\|_{\tilde{L}^1([T_1,t], C^{s-1})}\\
=C_F \max(\nu, \nu') \sup_{q\geq -1} 2^{q(s-1)} \int_{T_1}^t \|\D_q\bigg( \n R(\Ome, \GL \Xtol) +2\n R(^t\n \Xtol, M_{\nu, \nu', F} \n \Ome)\bigg)\|_{L^\infty} d\tau\\
\leq C_F \max(\nu, \nu') \sup_{q\geq -1} 2^{qs} \int_{T_1}^t \sum_{l\geq q-N_0 \atop l\geq -1} \sum_{\aa=-1}^1 \bigg(\|\D_l \Ome\|_{L^\infty} \|\D_{l+\aa} \GL \Xtol\|_{L^\infty}\\
+\|\D_{l+\aa}\n \Xtl\|_{L^\infty} \|\D_l M_{\nu, \nu', F} \n \Ome\|_{L^\infty}\bigg) d\tau\\
\leq C_F \max(\nu, \nu') \sup_{q\geq -1} 2^{qs} \sum_{l\geq q-N_0 \atop l\geq -1} \int_{T_1}^t 2^{2l} \|\D_l \Ome(\tau)\|_{L^\infty} 2^{-ls} \|\Xtol\|_{C^s} d\tau\\
\leq C_F \max(\nu, \nu') \left(\sum_{l\geq q-N_0} 2^{(q-l)s}\right) \|\Xtol\|_{L^\infty([T_1,t],C^s)} \|\Ome\|_{\tilde{L}^1([T_1,t], C_*^2)}.
\end{multline}
Similarly, for $F_\B$, thanks to Proposition \ref{propM1} (see Appendix) we can write that for all functions $f,g$ and all $l,l'\geq -1$:
$$
\|M(\D_l f, \D_{l'}g)\|_{L^\infty} \leq C_F 2^\frac{l+l'}{2} \|\D_l f\|_{L^\infty} \|\D_{l'} g\|_{L^\infty},
$$
which implies that (recall that $\La$ is a pseudo-differential operator of order $1$)
\begin{multline}
\|F_\B\|_{\tilde{L}^1([T_1,t], C^{s-1})}\\
\leq C_F |\nu-\nu'| \sup_{q\geq -1} 2^{qs} \sum_{l\geq q-N_0} \sum_{\aa=-1}^1 \int_{T_1}^t 2^{2l} \|\D_l \Xtol\|_{L^\infty} \|\D_{l+\aa}\Ome\|_{L^\infty} d\tau,
\end{multline}
wich gives the desired estimate. We now turn to the second term. Thanks to the the fact that for any functions $f,g$, the product $S_{q-1} f\cdot \D_q g$ has its frequencies localized in an annulus of size $2^q$ (see \eqref{supports}),
\begin{multline}
\|G_\A\|_{C^{s-3}} \leq \sup_{q\geq -1} 2^{q(s-2)}\bigg( \|S_{q-1}\Ome\|_{L^\infty} \|\D_l \GL \Xtl\|_{L^\infty}  +\|S_{q-1} \GL \Xtl\|_{L^\infty}  \|\D_q\Ome\|_{L^\infty}\\
+ 2\|S_{q-1}\n \Xtl\|_{L^\infty} \|M_{\nu, \nu', F} \n \D_q \Ome\|_{L^\infty}  +2\|S_{q-1} M_{\nu, \nu', F} \n \Ome\|_{L^\infty}  \|\n \D_q \Xtl\|_{L^\infty} \bigg)\\
\leq C_F \max(\nu, \nu') \sup_{q\geq -2} 2^{q(s-2)} \|\Ome\|_{L^\infty} \bigg( 2^{q} \|S_{q-1} \n\Xtl\|_{L^\infty} +2^{2q} \|\D_q \Xtl\|_{L^\infty}\bigg).
\end{multline}
The only difficulty is here that if we put all derivatives out of the norm $\|S_{q-1} \Xtl\|_{L^\infty}$ we cannot recover the $C^s$-norm of $\Xtl$. To overcome this problem we simply use \eqref{BesovS}, as $s-1<0$:
\begin{equation}
\|\n \Xtl\|_{C^{s-1}} \sim \sup_{q\geq -1} 2^{q(s-1)} \|\n S_q \Xtl\|_{L^\infty},
\label{trickSq}
\end{equation}
which implies
$$
\|G_\A\|_{C^{s-3}} \leq C_F \max(\nu, \nu') \sup_{q\geq -2} 2^{q(s-2)}  2^{2q} \|\Ome\|_{L^\infty} 2^{-qs} \|\Xtl\|_{C^s}.
$$
Proving the same for $G_\B$ will be more tricky for the same reason. More precisely, with the same method as before, some terms will involve $S_{q}$ and not enough derivatives for us to be able to use \eqref{trickSq}. To overcome this difficulty we will have to use Proposition \ref{propM2} (see Appendix). These estimates for the bilinear operator $M$ will enable us to move some derivatives where they will be needed to apply \eqref{trickSq}. With the same simplification as before we write,
\begin{multline}
G_\B=-(\nu-\nu')F^2(1-F^2)\n \sum_{q\geq -1}\bigg[ S_{q-1}\La\Xtl \cdot \D_q\La\Ome +S_{q-1} \La \Ome \cdot \D_q \La\Xtl\\
+M(S_{q-1} \La \Ome, \D_q \Xtli)+M(S_{q-1}\Xtl, \D_q\La \Ome) +\La\bigg(S_{q-1}\La\Xtl \cdot \D_q\Ome\\
+S_{q-1}\Ome \cdot \D_q\La\Xtl+M\big(S_{q-1}\Xtl, \D_q\Ome\big) +M\big(S_{q-1}\Ome, \D_q\Xtl\big) \bigg)\bigg].
\end{multline}
As for $G_\A$, using that for all $f,g$ we have (thanks again to Proposition \ref{propM1}),
$$
\|M(S_l f, \D_{l'}g)\|_{L^\infty} \leq C_F 2^\frac{l+l'}{2} \|S_l f\|_{L^\infty} \|\D_{l'} g\|_{L^\infty},
$$
allows us to get that:
\begin{multline}
\|G_\B\|_{C^{s-3}}\\
\leq C_F |\nu-\nu'|\sup_{q\geq -1} 2^{q(s-2)}\bigg(2^{q} \|S_{q-1} \La\Xtl\|_{L^\infty} \|\D_q \Ome\|_{L^\infty} +2^{2q} \|S_{q-1} \Ome\|_{L^\infty} \|\D_q\Xtl\|_{L^\infty}\\
+\|M(S_{q-1}\Xtl, \D_q\La \Ome)\|_{L^\infty} + \|\La M\big(S_{q-1}\Xtl, \D_q\Ome\big)\|_{L^\infty} \bigg).
\label{estGB1}
\end{multline}
In this expression, the first two terms will easily be estimated as explained before by
\begin{equation}
2^{-q(s-2)}\|\Xtl\|_{C^s} \|\Ome\|_{L^\infty},
\label{estGB2}
\end{equation}
and we will estimate the other terms thanks to Proposition \ref{propM2} as it will enable us to make some derivatives pound on $S_{q-1}\Xtl$, and use \eqref{trickSq}. As $s\in]0,1[$ let $\sigma>0$ such that $s+\sigma \in]0,1[$. Then thanks to \eqref{estM2}, we can write that
\begin{multline}
\|M(S_{q-1}\Xtl, \D_q\La \Ome)\|_{L^\infty} +\|\La M\big(S_{q-1}\Xtl, \D_q\Ome\big)\|_{L^\infty}\\
\leq C_F 2^q \|S_{q-1}\Xtl\|_{\dot{B}_{\infty, \infty}^{s+\sigma}} \|\D_q \Ome\|_{\dot{B}_{\infty, 1}^{1-(s+\sigma)}}.
\end{multline}
We emphasize that we deal here with homogeneous Besov norms and our results involve inhomogeneous norms. Then, thanks to the frequency localization we will use the following estimates: for any any function $u$, and any $\aa>0$:
$$
\begin{cases}
\vspace{1mm}
\|u\|_{\dot{B}_{\infty, \infty}^\aa} =\||D|^\aa u\|_{\dot{B}_{\infty, \infty}^0} \leq \||D|^\aa u\|_{L^\infty},\\
\|u\|_{\dot{B}_{\infty, 1}^\aa} \leq C\sqrt{\|u\|_{L^\infty} \||D|^{2\aa} u\|_{L^\infty}},
\end{cases}
$$
and deduce that
\begin{multline}
\|M(S_{q-1}\Xtl, \D_q\La \Ome)\|_{L^\infty} +\|\La M\big(S_{q-1}\Xtl, \D_q\Ome\big)\|_{L^\infty}\\
\leq C_F 2^q \|S_{q-1}|D|^{s+\sigma}\Xtl\|_{L^\infty} 2^{q\big(1-(s+\sigma)\big)}\|\D_q \Ome\|_{L^\infty}.
\end{multline}
Now we have enough derivatives pounding on $S_{q-1}$, indeed:
$$
\Xtl \in C^s \Longleftrightarrow |D|^{s+\sigma}\Xtl \in C^{-\sigma},
$$
so that we have the estimate:
$$
\|S_{q-1}|D|^{s+\sigma}\Xtl\|_{L^\infty} \leq C 2^{q\sigma} \||D|^{s+\sigma} \Xtl\|_{C^{-\sigma}} \leq C 2^{q\sigma} \|\Xtl\|_{C^s},
$$
which implies that:
\begin{multline}
\|M(S_{q-1}\Xtl, \D_q\La \Ome)\|_{L^\infty} +\|\La M\big(S_{q-1}\Xtl, \D_q\Ome\big)\|_{L^\infty}\\
\leq C_F 2^{q(2-s)} \|\Xtl\|_{C^s} \|\Ome\|_{L^\infty}.
\end{multline}
Plugging this estimates and \eqref{estGB2} into \eqref{estGB1} concludes the proof of the proposition. $\blacksquare$

\subsubsection{Study of $\Xtl(x, D)\Ome$ : estimates}

In this section we will recollect the previous results in order to bound the $C^{s-1}$-norm of $\Xtl(x, D)\Ome$. For more simplicity, we will denote it $\Xdo$. First as it satisfies system \eqref{transportregtg}, and as $s\in ]0,1[$, we can use Theorem \ref{thCs} from the appendix (see \cite{FCestimLp} for more details): as we have \eqref{estimvee} and $\div v_\ee=0$, there exists constants $C_F>0$ and $C_s>0$ such that if $T_2-T_1>0$ is so small that:
\begin{enumerate}
\item $C_F(T_2-T_1)^{\frac{1}{4}} \leq \no^3$,
\item $e^{C\int_{T_1}^{T_2} \|\n \Ue (\tau)\|_{L^\infty}}-1 \leq \frac{1}{C_F M_{visc}}$,
\item $T_2-T_1 +\int_{T_1}^{T_2} \|\n \Ue\|_{L^\infty} d\tau \leq C_{s,\no}$
\end{enumerate}
\begin{rem}
\sl{These conditions on $T_1,T_2$ are implied by what follows: there exists a constant $\cu =\cu(F,s,\nu, \nu')>0$ such that $T_1<T_2$ satisfy:
\begin{equation}
T_2-T_1 +\int_{T_1}^{T_2} \|\n \Ue\|_{L^\infty} d\tau \leq \cu.
\label{condT1T2}
\end{equation}
}
\end{rem}
If this condition is satisfied, then from Theorem \ref{thCs} there exists a constant $C_{\no,F}>0$ such that for all $t\in [T_1, T_2]$ satisfying condition \eqref{condT1T2},
\begin{multline}
\|\Xdo\|_{\tilde{L}^\infty([T_1,t], C^{s-1})}\\
\leq C_{\no, F} \left(\|\Xdo(T_1)\|_{C^{s-1}} +\|\Fe\|_{\tilde{L}^1([T_1,t],C^{s-1})} +\frac{1}{\no}\|\Ge\|_{\tilde{L}^\infty ([T_1,t],C^{s-3})} \right),
\label{estimaprioriOm}
\end{multline}
where, from system \eqref{transportregtg}, we have:
\begin{equation}
\begin{cases}
\Fe =(\nu-\nu')F \Xtl(x,D) \big(\d_3 \D \theta_\ee\big) +\Xtl(x,D) q_\ee +F_\A +F_\B,\\
\Ge =G_\A +G_\B.
\end{cases}
\end{equation}
Thanks to the previous section (see \ref{estimFG}), we have:
\begin{multline}
\|F_\A\|_{\tilde{L}^1([T_1,t], C^{s-1})} +\|F_\B\|_{\tilde{L}^1([T_1,t], C^{s-1})} +\frac{1}{\no}\|\Ge\|_{\tilde{L}^\infty ([T_1,t],C^{s-3})} \\
\leq C_F M_{visc} e^{T_1 D_F} \|\Xtol\|_{L^\infty([T_1,t],C^s)}.
\end{multline}
Next, there is no change for the following term, as in \cite{FC3} thanks to Besov product laws:
\begin{multline}
\|X^\ee(x,D) q_\ee\|_{\tilde{L}^1([T_1,t],C^{s-1})} \leq \int_{T_1}^t \|X^\ee\otimes q_\ee\|_{C^s} d\tau \leq \int_{T_1}^t \bigg(\|X^\ee\|_{C^s} \|\n \Ue\|_{L^\infty} \|\n \Uosc\|_{L^\infty}\\
+\|X^\ee\|_{L^\infty} (\|\n \Ue\|_{C^s} \|\n \Uosc\|_{L^\infty} +\|\n \Ue\|_{L^\infty} \|\n \Uosc\|_{C^s} \bigg)d\tau,
\end{multline}
then using \eqref{estimUoscTgamma}, \eqref{estapriori} and under condition \eqref{condmM}, we easily obtain:
\begin{equation}
\|X^\ee(x,D) q_\ee\|_{\tilde{L}^1([T_1,t],C^{s-1})} \leq C_F\big(1+\|\Xtol\|_{L^\infty([T_1,t], C^s)}\big).
\end{equation}
The additional term is also estimated thanks to product laws:
\begin{multline}
|\nu-\nu'|F \|\Xtl(x,D) \d_3 \D \theta_{\ee,osc}\|_{\tilde{L}^1([T_1,t],C^{s-1})} \leq C_F |\nu-\nu'|\int_{T_1}^t \|X^\ee. \d_3 \D \theta_{\ee, osc}\|_{C^s} d\tau\\
\leq C_F\big(1+\|\Xtol\|_{L^\infty([T_1,t], C^s)}\big).
\end{multline}
Substituting all these estimates in \eqref{estimaprioriOm} gives that, under conditions \eqref{condT1T2} and \eqref{condmM}, if $0<\ee\leq \ee_0$, for all $t\in [T_1, T_2]$,
\begin{multline}
\|\Xdo\|_{\tilde{L}^\infty([T_1,t], C^{s-1})}\\
\leq C_{\no, F} \left(\|\Xdo(T_1)\|_{C^{s-1}} +1 +e^{T_1 D_F} \sup_{\tau \in [T_1,t]}\|\Xtol\|_{C^s} \right),
\label{estimXom}
\end{multline}
In \cite{FC3} after having obtained the analoguous of \eqref{estimXe2} and \eqref{estimXom} we could easily obtain a globalization for large times of the result. In our case, we have to be extremely careful: due to \eqref{estimOmLinf}, the $L^\infty$ estimates of the potential vorticity involves the coefficient $e^{T_1 D_F}$ that may be very large. For this reason the only way we can absorb the last term of \eqref{estimXom} is to consider:
$$
Z(t) \overset{def}{=} \|\Xdo\|_{\tilde{L}^\infty([T_1,t], C^{s-1})} +2 C_{\no, F} e^{T_1 D_F} \|\Xtol\|_{L^\infty( [T_1,t], C^s)}.
$$
Then combining \eqref{estimXe2} and \eqref{estimXom} we get that for all $t\in[T_1,T_2]$,
$$
Z(t) \leq C_{\no, F} \bigg(Z(T_1) +1 + \int_{T_1}^t \big(e^{T_1 D_F} +\|\n \Ue\|_{L^\infty} \big) Z(t') dt'\bigg) +\frac{1}{2} Z(t),
$$
Simplifying and thanks to the Gronwall estimate, we end up with (thanks to condition \eqref{condT1T2}):
\begin{multline}
Z(t) \leq C_{\no, F} \big(Z(T_1) +1\big) e^{\displaystyle{C_{\no, F}\int_{T_1}^t \big(e^{T_1 D_F} +\|\n \Ue\|_{L^\infty} \big) dt'}}\\
\leq C_{\no, F} \big(Z(T_1) +1\big) e^{\displaystyle{C_{\no, F} (e^{T_1 D_F}+ \int_{T_1}^t \|\n \Ue\|_{L^\infty} \big) dt'}}.
\end{multline}
The globalisation argument is classical except that we have to be careful with $e^{T_1 D_F}$: for all $t\in [0,\min(T_\ee, T_\ee^\gamma)$, we subdivide $[0,t]$ into $0=T_0 <T_1<...<T_N=t$ such that for every $i\in \{0,N-1\}$,
\begin{equation}
T_{i+1}-T_i +\int_{T_i}^{T_{i+1}} \|\n \Ue\|_{L^\infty} dt' \sim \cu.
\label{condTi2}
\end{equation}
The previous arguments give us that for all $i\in \{0,N-1\}$,
\begin{multline}
Y_{i+1} \overset{def}{=}\|\Xdo\|_{\tilde{L}^\infty([T_i,T_{i+1}], C^{s-1})} +2 C_{\no, F} e^{T_i D_F} \|\Xtol\|_{L^\infty( [T_i,T_{i+1}], C^s)}\\
\leq C_{\no, F} \big(Y_i +1\big) e^{\displaystyle{C_{\no, F} (e^{T_i D_F}+ \int_{T_i}^{T_{i+1}} \|\n \Ue\|_{L^\infty} \big) dt'}},
\end{multline}
and by induction (using that $t=T_N$) we get that:
\begin{multline}
Y(t) \overset{def}{=} \|\Xdo\|_{\tilde{L}^\infty([0,t], C^{s-1})} +\|\Xtol\|_{L^\infty( [0,t], C^s)}\\
\leq \|\Xdo\|_{\tilde{L}^\infty([0,t], C^{s-1})} +2 C_{\no, F} e^{T_{N-1} D_F} \|\Xtol\|_{L^\infty( [0,t], C^s)}\\
\leq C_{\no, F}^N (1+Y(0)) e^{\displaystyle{C_{\no, F} \big( e^{T_0 D_F} +...+ e^{T_{N-1} D_F}\big)}} e^{\displaystyle{\int_0^t \|\n \Ue\|_{L^\infty}dt'}}\\
\leq C_0 C_{\no, F}^N e^{\displaystyle{C_{\no, F} N e^{t D_F}}} e^{\displaystyle{\int_0^t \|\n \Ue\|_{L^\infty} dt'}}.
\end{multline}
Finally, summing \eqref{condTi2} for $i\in \{0,N-1\}$ we obtain that:
$$
N \sim \frac{1}{\cu} \left(t+ \int_0^t \|\n \Ue\|_{L^\infty} dt' \right).
$$
Substituting this in the previous result and estimating the oscillating part as before finally gives that:
\begin{multline}
Y(t) \leq C_0 e^{\displaystyle{N\left(\ln C_{\no, F} + C_{\no, F} e^{t D_F}\right)}} e^{\displaystyle{\int_0^t \|\n \Ue\|_{L^\infty} dt'}}\\
\leq C_0 e^{\displaystyle{2C_{\no, F} e^{t D_F}\left(t+ \int_0^t \|\n \Ue\|_{L^\infty} dt' \right)}} e^{\displaystyle{\int_0^t \|\n \Uqg\|_{L^\infty} dt'}}\\
\leq C_0 e^{\displaystyle{3C_{\no, F} e^{t D_F}\left(t+ \int_0^t \|\n \Ue\|_{L^\infty} dt' \right)}}.
\end{multline}
We can estimate the oscillating part: thanks to \eqref{estimUoscsimple}
$$
\int_0^t \|\n \Uosc\|_{L^\infty} dt' \leq t^\frac{3}{4} \|\n \Uosc\|_{L_t^4 L^\infty} \leq C_F (T_\ee^\gamma)^\frac{5}{4} \ee^\frac{M}{2} \leq 1,
$$
if $0< \ee \leq \ee_0$, with $\ee_0$ small enough, we end up with:
\begin{multline}
Y(t) \overset{def}{=} \|\Xdo\|_{\tilde{L}^\infty([0,t], C^{s-1})} +\|\Xtol\|_{L^\infty( [0,t], C^s)} \\
\leq C_0 e^{\displaystyle{3C_{\no, F} e^{t D_F}\left(t+1+ \int_0^t \|\n \Uqg\|_{L^\infty} dt' \right)}}.
\label{estimNUqg}
\end{multline}

\subsection{End of the proof}

We are now able to conclude the argument started with \eqref{bootstrap}. First, thanks to the logarithmic estimate \eqref{estimlog2}
\begin{equation}
\|\n \Uqg\|_{L^{\infty}} \leq C\Bigg(\|\Ome\|_{L^2}+ \|\Ome\|_{L^{\infty}} \log \big(e+1+\frac{\|[\Xtl]^{-1}\|_{L^{\infty}} +Y(t) \big)}{\|\Ome\|_{L^{\infty}}} \big) \Bigg).
\end{equation}
Using that for any $A,B>0$ the function $x\mapsto x\ln (A+\frac{B}{x})$ is increasing on $]0,\infty[$, and substituting \eqref{estimOmegaL2}, \eqref{estimOmLinf}, \eqref{estimXe} and \eqref{estimNUqg} into the previous estimate, there exists $\ee_0,C_F>0$ and $M$ such that if condition \eqref{condmM} is satisfied, then for all $t\leq \min(T_\ee, T_\ee^\gamma)$,
\begin{multline}
\|\n \Uqg\|_{L^{\infty}} \leq C_F\Bigg(1+ e^{tD_F} \log \big(e+1+C_0 e^{\displaystyle{3C_{\no, F} e^{t D_F}\big(t+1+ \int_0^t \|\n \Uqg\|_{L^\infty} dt' \big)}}\big) \Bigg)\\
\leq C_F\Bigg(1+ e^{tD_F} \log C_0 +3C_{\no, F} e^{2t D_F}\big(t+1+ \int_0^t \|\n \Uqg\|_{L^\infty} dt' \big)\Bigg)\\
\leq C_{\no, F} te^{2tD_F} +3 C_{\no, F} e^{2tD_F} \int_0^t \|\n \Uqg\|_{L^\infty} dt'.
\end{multline}
The Gronwall estimates then implies that for all $t\leq \min(T_\ee, T_\ee^\gamma)$,
\begin{equation}
\|\n \Uqg\|_{L^{\infty}} \leq C_{\no, F} te^{2tD_F} e^{\displaystyle{3 C_{\no, F} te^{2tD_F}}} \leq e^{\displaystyle{4 C_{\no, F} te^{2tD_F}}}.
\end{equation}
\begin{rem}
\sl{We can afford such a rough estimates as $x \leq e^x$, because the precision is forced by the term of size $e^{e^t}$.
}
\end{rem}
Following \eqref{bootstrap}, for all $t\leq \min(T_\ee, T_\ee^\gamma)$, we have $V_\ee(t) \leq K_\ee=-2\ln \ee$ and thanks to \eqref{estimUoscsimple}:
\begin{multline}
V_\ee(t) =\int_0^t \|\n \Ue(t')\|_{L^\infty} dt' \leq (T_\ee^\gamma)^\frac{5}{4} \ee^\frac{M}{2} +T_\ee^\gamma e^{\displaystyle{4 C_{\no, F} T_\ee^\gamma e^{2D_F T_\ee^\gamma}}} \leq 2T_\ee^\gamma e^{\displaystyle{4 C_{\no, F} T_\ee^\gamma e^{2D_F T_\ee^\gamma}}} \\
\leq 2\gamma \ln(\ln |\ln \ee|)  e^{\displaystyle{4 C_{\no, F} \gamma \ln(\ln |\ln \ee|) (\ln |\ln \ee|)^{2D_F \gamma}}}
\end{multline}
If $\gamma$ is so small that $2D_F \gamma <1$ and $4 C_{\no, F} \gamma <\frac{1}{2}$, as there exists a constant $C$ such that for all $x\geq e$,
$$
\ln(x) e^{\displaystyle{4 C_{\no, F} \gamma \ln (x) x^{2D_F \gamma}} }\leq Ce^{\frac{x}{2}},
$$
we end up with:
\begin{equation}
V_\ee(t) \leq 2 \gamma C|\ln \ee|^\frac{1}{2} \leq \frac{K_\ee}{2},
\end{equation}
if $\ee$ is small enough, which is the bound announced in \eqref{bootestim} that allows us to obtain that $T_\ee^* >T_\ee^\gamma$. The last part of the proof is done exactly as in \cite{FC3}: indeed the non-local operator is not a problem for $L^2$-estimates. This concludes the proof of the theorem. $\blacksquare$

\section{Appendix}

The first part is devoted to a quick presentation of the Littlewood-Paley theory. In the second section we briefly recall general definitions for vortex patches and the last section provides new properties for the operator $\G$ and recalls the a priori estimates from \cite{FCestimLp}.

\subsection{Littlewood-Paley theory}

In this section, $C^s$ is the usual H\"older space, which can also be defined through the Littlewood-Paley theory if $s\notin \N$ (we refer to \cite{Chemin1, Dbook} for a complete presentation of the theory)~:
$$
C^s=\{u \in \cS'(\R^3), \quad \|u\|_{C^s} \overset{\mbox{def}}{=} \sup_{q \geq -1} 2^{qs} \|\Delta_q u\|_{L^{\infty}} < \infty\},  
$$
where $\Delta_q$ is the classical dyadic frequency localization operator defined as follows~: consider a smooth radial function $\chi$ supported in the ball $B(0, \frac{4}{3})$, equal to $1$ in a neighborhood of $B(0, \frac{3}{4})$ and such that $r\mapsto \chi(r.e_1)$ is nonincreasing over $\R_+$. If we define $\varphi(\xi)=\chi(\frac{\xi}{2})-\chi(\xi)$, $\varphi$ is supported respectively in the annulus $\cC(0,\frac{3}{4}, \frac{8}{3})$ (equal to $1$ in a sub-annulus), and satisfiy that for all $\xi\in \R^3$,
$$
\chi(\xi)+ \sum_{q\geq 0} \varphi(2^{-q} \xi)=1 \quad \mbox{and if }\xi \neq 0, \quad \sum_{q\in\Z} \varphi(2^{-q} \xi)=1.
$$
Then for all tempered ditribution we define :
\begin{itemize}
\item $\Delta_{-1}= \mathcal{F}^{-1}\big(\chi(\xi) \hat{u}(\xi)\big)$ and $\forall q\leq -2, \quad \Delta_q=0$,
\item $\displaystyle{\forall q\geq 0, \quad \Delta_q =\mathcal{F}^{-1}\big(\varphi(2^{-q}\xi) \hat{u}(\xi)\big)}$ and $\displaystyle{S_q u =\sum_{p<q-1} \Delta_p u= \chi(2^{-q}D) u}$,
\item $\displaystyle{\forall q \in \mathbb{Z}, \quad \dot{\Delta}_q =\mathcal{F}^{-1}\big(\varphi(2^{-q}\xi) \hat{u}(\xi)\big)}$ and $\displaystyle{\dot{S}_q u  =\sum_{p<q-1} \dot{\Delta}_p u= \chi(2^{-q}D) u}$.
\end{itemize}
H\"older spaces are particular cases of inhomogeneous Besov spaces: $C^s=B_{\infty, \infty}^s$, where
$$
B_{p,r}^s=\{u \in \cS'(\R^3), \quad \|u\|_{B_{p,r}^s} \overset{\mbox{def}}{=} \|\big(2^{qs} \|\Delta_q u\|_{L^p} \big)_{q \geq -1}\|_{\ell^r} < \infty\}.
$$
The homogeneous Besov spaces are defined as follows:
$$
\dot{B}_{p,r}^s=\{u \in \cS'(\R^3), \quad \|u\|_{\dot{B}_{p,r}^s} \overset{\mbox{def}}{=} \|\big(2^{qs} \|\ddq u\|_{L^p} \big)_{q \in \Z}\|_{\ell^r} < \infty\}.
$$
When the regularity index $s$ is negative, another way to express the Besov norm involves the operator $S_q$ instead of $\D_q$ (for more details we refer to \cite{Dbook} Proposition $2.76$ and to $2.31$ for the homogeneous case):
\begin{prop}
\sl{There exists a constant $C>0$ such that for all $s<0$, $p,r\in[1,\infty]$ and $u$, then $u\in B_{p,r}^s$ if and only if
$$
\big(2^{qs} \|S_q u\|_{L^p}\big)_{q \geq -1} \in \ell^r.
$$
Moreover, we have
\begin{equation}
\|\big(2^{qs} \|S_q u\|_{L^p}\big)_{q}\|_ {\ell^r} \sim \|u\|_{B_{p,r}^s}.
\label{BesovS}
\end{equation}
}
\end{prop}
\begin{rem}
\sl{Due to the supports, we easily obtain that
\begin{equation}
\D_j \D_l =\dot{\D}_j \dot{\D}_l =0 \quad \mbox{if } |j-l|\geq 2.
\label{supports}
\end{equation}
\begin{itemize}
\item For any functions $f,g$, and any $\aa\in\{-1,0,1\}$, the product $\D_q f. \D_{q+\aa} g$ has its frequencies in a ball of size $2^q$.
\item For any functions $f,g$, the product $S_{q-1} f. \D_q g$ has its frequencies in an annulus of size $2^q$.
\end{itemize}
}
\end{rem}
We will end this section with the Bony decomposition, which comes from the fact that for all distributions $u,v$, we can write the product as follows:
$$
uv=(\sum_{q\geq -1} \D_q u)(\sum_{l\geq -1} \D_l v).
$$
In fact, a more efficient way to write this product is the following Bony decomposition, where we basically set three parts according to the fact that the frequency $q$ of $u$ is smaller, comparable or bigger than the frequency $l$ of $v$:
\begin{equation}
uv= T_u v+ T_v u+ R(u, v),
\label{Bonydecp}
\end{equation}
where
\begin{itemize}
\item $T$ is the paraproduct~: $T_u v= \sum_{p\leq q-2} \Delta_p u \Delta_q v= \sum_{q} S_{q-1} u \Delta_q v$,
\item $R$ is the remainder~: $R(u, v)= \sum_{|p-q| \leq 1} \Delta_p u \Delta_q v$.
\end{itemize}
A similar decomposition can be defined with the homogeneous Littlewood-Paley operators.
\label{appendiceLP}

\subsection{Vortex patches}

We refer to \cite{Chemin1} for a full description of the persistence of the vortex patches structure in the case of the Euler system, to \cite{Dpoches} and \cite{TH1} for the case of the Navier-Stokes system, and to \cite{Dutrifoy1} for the case of the inviscid Primitive Equations ($\nu=\nu'=0$). In the present paper, we take here the same definitions of vortex patches and tangential regularity as in \cite{Dutrifoy1}~: a potential vortex patch will be defined with respect to the scalar potential vorticity instead of the vorticity (rotational of the velocity). The potential vorticity is called a vortex patch if it is the characteristic function of a regular open set~:
\begin{defi}
\sl{
We say that $\Omega_0$ is a vortex patch of class $C^s$ if, for some $s\in ]0,1[$,
$$
\Om_0= \Om_{0, i} \textbf{1}_D +\Om_{0, e} \textbf{1}_{\R^3-D},
$$
where $\Om_{0, i} \in C^s(\overline{D})$, $\Om_{0, e} \in C^s(\R^3-D)$ and $D$ is an open bounded domain of class $C^{s+1}$.
}
\end{defi}
In our results we will use estimates involving the tangential regularity with respect to a set $X$ of vectorfields :
\begin{defi}
\sl{If $X=(X_{\lambda})_{\lambda=1,...,N}$ is a finite family of vectorfields we will say that this family is admissible if and only if ($\times$ is the usual vector product in $\mathbb{R}^3$) :
$$
[X]^{-1}\overset{\mbox{def}}{=} \Big(\frac{2}{N(N-1)} \sum_{\lambda< \lambda'} |X_{\lambda} \times X_{\lambda'}|^2 \Big)^{-\frac{1}{4}} <
\infty.
$$
If $s\in ]0,1[$ and $X$ is an admissible family of vectorfields $C^s$ we define the space :
$$
C^s(X)=\{w \in L^{\infty} \mbox{ such that} \quad\ X_{\lambda}(x,
D)w \overset{\mbox{def}}{=} \mbox{div} (w \otimes X_{\lambda})\in C^{s-1} \}
$$
and as corresponding norm we take~:
\begin{equation}
\|w\|_{C^s(X)}\overset{\mbox{def}}{=}\|w\|_{L^{\infty}}+\|[X]^{-1}\|_{L^{\infty}} +\sum_{\lambda=1}^N \big( \|X_{\lambda}\|_{C^s}+ \|X_{\lambda}(x, D)w\|_{C^{s-1}} \big).
\end{equation}
}
\end{defi}
\begin{rem}
\sl{We took here the same definition as in \cite{Dutrifoy1} for $X_{\lambda}(x, D)w$ which is a slightly simplified formulation of the definitions from \cite{Chemin1}, or \cite{TH1}.}
\end{rem}

\subsection{Definition and additional properties for the non-local operator $\G$}

\subsubsection{Product estimates}

We refer to \cite{FCestimLp} for a study of operator $\G$ where we give various formulations. This operator is defined as
$$
\G= \D \DF^{-1} (\nu \d_1^2 +\nu \d_2^2+ \nu' F^2 \d_3^2).
$$
First we decompose $\G$ into its local and non-local parts: 
$$
\G =\GL +(\nu-\nu') F^2(1-F^2) \La^2,
$$
where we denote:
\begin{equation}
\begin{cases}
\GL= \nu \d_1^2 +\nu \d_2^2+ \left((1-F^2)\nu +F^2 \nu'\right) \d_3^2,\\
\La= \d_3^2 (-\DF)^{-\frac{1}{2}}.
\end{cases}
\label{decpgamma}
\end{equation}
We also refer to \cite{FCestimLp} for the following expressions of $\La$ as singular or convergent integrals, directly related to the alternative expression of homogeneous Besov norms involving finite differences:
\begin{multline}
\La f (x) =\lim_{\ee \rightarrow 0} \int_{|y|\geq \ee}  K(y)\big(f(x-y)-f(x)\big) dy\\
=\frac{1}{2}\int_{\R^3}  K(y)\bigg(f(x-y)+f(x+y)-2f(x)\bigg) dy,
\end{multline}
where the kernel $K$ is defined for all $y\in \R^3$ by ($C$ is a universal constant):
\begin{equation}
K(y)= -\frac{2C}{F^3}\frac{y_1^2+ y_2^2 -\frac{3}{F^2}y_3^2}{\big(y_1^2+ y_2^2 +\frac{1}{F^2}y_3^2\big)^3}.
\label{Kexpr}
\end{equation}
The main feature of \cite{FCestimLp} was an estimate of the commutator of $\La$ with a lagrangian change of variable (crucial to obtain the apriori estimates), but we also obtained the following result that allows to consider commutators like $\G(fg)-f\G g$:
\begin{prop}
\sl{(\cite{FCestimLp} Section $3.3.2$) For any smooth functions $f,g$ we can write:
$$
\La (fg)= f\La g +g\La f +M(f,g),
$$
where the bilinear operator $M$ is defined for all $x\in \R^3$ by:
\begin{equation}
M(f,g)(x) =\int_{\R^3} K(y) \big(f(x-y)-f(x)\big) \big(g(x-y)-g(x)\big) dy.
\label{defM}
\end{equation}
Moreover there exists a constant $C_F$ such that for all $f,g$:
\begin{equation}
\|M(f,g)\|_{L^p} \leq C_F \sqrt{\|f\|_{L^p} \|\n f\|_{L^p} \|g\|_{L^\infty} \|\n g\|_{L^\infty}},
\end{equation}
}
\label{propM1}
\end{prop}
In this article we will need more precise estimates where we can in particular, even if $M(f,g)=M(g,f)$, make the derivatives pound differently on $f$ or $g$, which is the object of the following result:
\begin{prop}
\sl{There exists a constant $C_F$ such that for all $f,g$ and all $p,p_1, P_2,r,\overline{r}\in[1,\infty]$ and $\eta$ satisfying:
$$
\begin{cases}
\vspace{1mm}
\displaystyle{\frac{1}{p}= \frac{1}{p_1} +\frac{1}{p_2}, \quad 1= \frac{1}{r} +\frac{1}{\overline{r}},}\\
\displaystyle{2-\eta-\frac{3}{r} \in]0,1[,}
\end{cases}
$$
then we have
\begin{equation}
\|M(f,g)\|_{L^p} \leq C_F \|f\|_{\displaystyle{\dot{B}_{p_1,r}^{2-\eta-\frac{3}{r}}}}\|g\|_{\displaystyle{\dot{B}_{p_2,\overline{r}}^{2+\eta-\frac{3}{\overline{r}}}}}.
\end{equation}
}
\label{propM2}
\end{prop}
\begin{rem}
\sl{In the present article we will use the previous proposition in the case $p=p_1=p_2=\infty$, $r=\infty, \overline{r}=1$ and for $\eta=2-(s+\gamma)$ where $s\in]0,1[$ and $\sigma>0$ is such that $s+\sigma<1$. In this case we end up with:
\begin{equation}
\|M(f,g)\|_{L^\infty} \leq C_F \|f\|_{\dot{B}_{\infty,\infty}^{s+\sigma}}\|g\|_{\dot{B}_{\infty,1}^{1-(s+\sigma)}}.
\label{estM2}
\end{equation}
}
\end{rem}
\textbf{Proof :} let $f,g$ be smooth functions. From the expression of the kernel, there exists a constant $C_F>0$ such that, thanks to the relations between the parameters:
\begin{multline}
\|M(f,g)\|_{L_x^p} \leq C_F \int_{\R^3} \frac{\|f(.-y)-f(.)\|_{L_x^{p_1}}}{|y|^{2-\eta}} \cdot \frac{\|g(.-y)-g(.)\|_{L_x^{p_2}}}{|y|^{2+\eta}} dy\\
\leq C_F \left(\int_{\R^3} \frac{\|f(.-y)-f(.)\|_{L_x^{p_1}}^r}{|y|^{(2-\eta)r}} dy\right)^\frac{1}{r} \left( \int_{\R^3} \frac{\|g(.-y)-g(.)\|_{L_x^{p_2}}^{\overline{r}}}{|y|^{(2+\eta)\overline{r}}} dy \right)^\frac{1}{\overline{r}}.
\end{multline}
If we denote $s_1=2-\eta-\frac{3}{r}$ and $s_2=2+\eta-\frac{3}{\overline{r}}$ then $(2-\eta)r= s_1 r+3$ and $(2+\eta)\overline{r} =s_2 \overline{r} +3$. Moreover $s_1+s_2=1$ so if $1-\frac{3}{r}<\eta<2-\frac{3}{r}$ both regularity exponents are in $]0,1[$ and we can use the following result
\begin{thm}\sl{(\cite{Dbook}, $2.36$)
 Let $s \in ]0,1[$ and $p,r\in [1,\infty]$. There exists a constant $C$ such that for any $u\in \cC_h'$,
$$
C^{-1} \|u\|_{\dot{B}_{p,r}^s}\leq \|\frac{\|\tau_{-y}u -u\|_{L^p}}{|y|^s}\|_{L^r (\R^d; \frac{dy}{|y|^d})} \leq C \|u\|_{\dot{B}_{p,r}^s}.
$$}
\end{thm}
This concludes the proof. $\blacksquare$

\subsubsection{A priori estimates}

For the confort of the reader we state here the a priori estimates obtained in \cite{FCestimLp}. We consider the following transport diffusion system:

\begin{equation}
\begin{cases}
\d_t u +v.\n u -\G u =\Fe,\\
u_{|t=0}=u_0
\end{cases}
\label{TDQG}
\end{equation}
Let us introduce $M_{visc}=\frac{\max(\nu, \nu')}{\min(\nu, \nu')}$.

\begin{prop}\sl{($L^p$-estimates) Assume that $u$ solves \eqref{TDQG} on $[0,T]$ with $u_0 \in L^p$ and that $\|v\|_{L_T^\infty L^6} \leq C'$ (for some constant $C'$) with $\div v=0$. Then there exists a constant $D$ (depending on $F$, $M_{visc}$ and $C'$) such that for all $t\in [0,T]$,
\begin{equation}
\|u\|_{L_t^\infty L^p} \leq D^t (\|u_0\|_{L^p}+\int_0^t \|\Fe(\tau)\|_{L^p} d\tau).
\end{equation}  
}
\label{appendicestimLp}
\end{prop}

\begin{thm}\sl{(Smoothing effect) Assume that $u$ solves \eqref{TDQG} on $[T_1, T_2]$ with $v$ satisfying $\div v=0$ and $\|v\|_{L^\infty({[T_1,T_2]}, L^6)} \leq C'$, $u(T_1)\in L^p$, $\Fe \in L_{loc}^1 L^p$ (for $p\in[1,\infty]$). There exist two constants $C$ and $C_F$ such that if $T_2-T_1>0$ is so small that:
\begin{enumerate}
\item $2CC'(T_2-T_1)^{\frac{1}{4}} \leq \no^\frac{3}{4}$,
\item $e^{C\int_{T_1}^{T_2} \|\nabla v (\tau)\|_{L^\infty}}-1 \leq \frac{1}{C_F M_{visc}}$.
\end{enumerate}
Then, for all $r\in[1,\infty]$, there exists a constant $C_{r,F}>0$ such that for all $t\in [T_1, T_2]$,
\begin{equation}
(\no r)^{\fr} \|u\|_{\tilde{L}^r([T_1,t], B_{p, \infty}^\frb)} \leq C_{r, F} \left(\|u(T_1)\|_{L^p} +\int_{T_1}^t \|\Fe(\tau)\|_{L^p} d\tau \right).
\end{equation}
}
\label{thC2}
\end{thm}

\begin{rem}
\sl{In the particular case $r=1, p=\infty$ we obtain that:
\begin{equation}
\no \|u\|_{\tilde{L}^1([T_1,t], C_*^2)} \leq C_F \left(\|u(T_1)\|_{L^\infty} +\int_{T_1}^t \|\Fe(\tau)\|_{L^\infty} d\tau \right).
\end{equation}
}
\end{rem}

\begin{thm}\sl{(a priori estimates)
Let $s\in ]-1,1[$. Assume that $u$ solves \eqref{TDQG} on $[T_1, T_2]$ with $v$ satisfying $\div v=0$ and $\|v\|_{L^\infty({[T_1,T_2]}, L^6)} \leq C'$, $u(T_1)\in B_{p,\infty}^s$. Assume in addition that the external force term can be decomposed into $\Fe+\Ge$, with $\Fe \in \tilde{L}^1([T_1,T_2], B_{p,\infty}^s)$ (for $p\in[1,\infty]$) and $\Ge \in \tilde{L}^\infty([T_1,T_2], B_{p,\infty}^{s+\frac{2}{r}-2})$ for $r\in[1,\infty]$ with $s+\frac{2}{r}\in]-1,1[$. There exist two constants $C_s$ and $C_F$ such that if $T_2-T_1>0$ is so small that:
\begin{enumerate}
\item $2CC'(T_2-T_1)^{\frac{1}{4}} \leq \no^\frac{3}{4}$,
\item $e^{C\int_{T_1}^{T_2} \|\nabla v (\tau)\|_{L^\infty}}-1 \leq \frac{1}{C_F M_{visc}}$,
\item $T_2-T_1 +\int_{T_1}^{T_2} \|\n v\|_{L^\infty} d\tau \leq C_{s, \no}$
\end{enumerate}
Then, there exists a constant $C_{\no,F}>0$ such that for all $r\in [1,\infty]$ with $s+\frac{2}{r}\in]-1,1[$ and  $t\in [T_1, T_2]$,
\begin{equation}
(\no r)^\frac{1}{r}\|u\|_{\tilde{L}^r([T_1,t], B_{p, \infty}^{s+\frac{2}{r}})} \leq C_{\no, F} \left(\|u(T_1)\|_{B_{p,\infty}^s} +\|\Fe\|_{\tilde{L}^1([T_1,t],B_{p,\infty}^s)} +\frac{1}{\no}\|\Ge\|_{\tilde{L}^\infty ([T_1,t],B_{p,\infty}^{s-2})} \right).
\end{equation}
}
\label{thCs}
\end{thm}

The author wishes to thank R. Danchin, I. Gallagher, and T. Hmidi for useful discussions.

\end{document}